\documentclass[reqno,oneside,11pt]{amsart}
\usepackage[french]{babel}
\usepackage{amscd,amssymb,epsfig,verbatim,xypic,fancyhdr,amsmath}
\usepackage[latin1]{inputenc}
\usepackage{graphicx}
\usepackage{psfrag}
\usepackage{yhmath}
\usepackage{frbib}

\title[Paysage systolique en caractéristique~-1]{Paysage systolique des surfaces hyperboliques compactes de caractéristique -1}
\author{Matthieu Raphaël Gendulphe}
\date{Le \today}

\theoremstyle{plain}
\newtheorem{thm}{Th\'eor\`eme}[section] 
\newtheorem{cor}[thm]{Corollaire}
\newtheorem{pro}[thm]{Proposition}
\newtheorem{lem}[thm]{Lemme}
\newtheorem{thm-defi}[thm]{Théorème - Définition}
\newtheorem{pro-defi}[thm]{Proposition - Définition}

\theoremstyle{definition}

\newtheorem{rem}{Remarque}[section]

\def\R{\mathbb{R}} 
 
\def\Z{\mathbb{Z}}

\def\t{\mathcal{T}} 
\def\T{\mathbb{T}}

\def\D{\rm{D}}
\def\H{\mathcal{H}}

\def\g{\gamma} 
\def\d{\delta}

\def\f{\varphi}

\def\s{\sigma}

\def\1{\bold{1}}

\def\Aut{\text{Aut}\,}

\xyoption{all}

\addtocounter{section}{0}              
\numberwithin{equation}{section}       

\begin{document}

\renewcommand{\contentsname}{Plan de l'article} 

\maketitle
\begin{abstract}
Nous déterminons des inégalités optimales pour la systole dans le cas des surfaces hyperboliques compactes de caractéristique -1.
Dans un premier temps nous étudions la géométrie de ces surfaces, nous
décrivons ensuite l'action des groupes modulaires sur les espaces de Teichmüller. Nous décomposons alors des domaines fondamentaux de ces actions en cellules adaptées à la systole, c'est-à-dire telles que tous les points d'une cellule ont mêmes géodésiques pour systoles.
Enfin, nous donnons tous les points critiques des fonctions systole et en particulier
leurs maxima. Nous nous intéressons aussi à d'autres invariants comme la 2-systole, la 3-systole...
\end{abstract}

\renewcommand{\abstractname}{Abstract} 
\tableofcontents{}

\subsection*{Introduction}
 La \emph{systole} d'une variété riemanienne compacte est définie comme la longueur minimale d'une courbe non contractile, par abus de langage on emploie le même mot pour les courbes réalisant cette longueur.\par
 La systole des surfaces de Riemann fut sujet d'études approfondies durant la décennie
 1993-2003. L'article de P.~Schmutz Schaller \cite{schmutz}, marqua une avancée importante dans la recherche de surfaces maximales pour la systole. Celui de C.~Bavard \cite{bavard1}, fournit quant-à lui un cadre théorique agréable, calqué sur celui de la théorie des réseaux.\par
 Nous nous inscrivons dans la continuité de ces travaux, néanmoins nous abordons ici
une direction jusqu'à présent délaissée~: celle des surfaces hyperboliques non orientables.
Nous nous intéressons au cas des surfaces compactes de caractéristique -1, c'est-à-dire
au plan projectif à deux bords, à la bouteille de Klein à un bord et à la somme connexe de
trois plans projectifs.
Pour chacune d'elles nous donnons une décomposition cellulaire d'un domaine fondamental de l'action du roupe modulaire, cette décomposition est adaptée à la systole car dans chaque cellule l'ensemble des géodésiques systoles ne varie pas. Il est alors facile de déterminer les points critiques de la systole sur l'espace
de Teichmüller correspondant. En particulier, nous en déduisons tous les maxima locaux et globaux de la systole, voici un extrait de nos principaux résultats~:

\begin{thm} Nous avons les inégalités optimales suivantes~:
\begin{itemize}
\item soient $X$ une surface hyperbolique fermée de caractéristique -1 et $s$ sa systole, alors $$\cosh(s)\leq\frac{5+\sqrt{17}}{2}~;$$
\item soient $Y$ une bouteille de Klein à un bord de longueur $b_1$ et $s$ sa systole, alors $$\cosh(s)\leq\cosh(\frac{b_1}{4})+1~;$$
\item soient $Z$ un plan projectif à deux bords de longueurs $b_1$, $b_2$ et $s$ sa systole, alors
$$\cosh(s)\leq\cosh(\frac{b_1}{2})+\cosh(\frac{b_2}{2})+1~;$$
\item soient $\T$ un tore à un bord de longueur $b_1$ et $s$ sa systole, alors
$$\cosh(\frac{s}{2})\leq\cosh(\frac{b_1}{6})+\frac{1}{2}.$$
\end{itemize}
\end{thm}

Dans le cas du tore à un bord l'inégalité était déjà connue (voir \cite{schmutz}).

\subsection*{Conventions} Dans toute la suite le terme \emph{surface$_{-1}$}
désignera une surface hyperbolique fermée de caractéristique -1.
Par \emph{géodésique} nous entendrons, sauf mention du contraire,
une géodésique fermée simple. Toute courbe fermée simple sera qualifiée d'\emph{orientable} si elle admet
un voisinage orientable, de \emph{non orientable} dans le cas contraire.
Nous noterons $S^{or}$ le revêtement des orientations d'une surface $S$, et $\g^{or}$ un relevé à $S^{or}$
d'une géodésique $\g$ de $S$.

 Les bords des surfaces hyperboliques seront toujours supposés totalement géodésiques.
L'\emph{auto-recollement} d'un bord consistera en l'identification
de ses points opposés. Lorsqu'il n'y aura pas d'ambiguité, nous conserverons la notation d'une géodésique après
opération de découpe ou d'auto-recollement. Les polygones évoqués tout au long du texte sont hyperboliques.
Nous confondrons les noms des géodésiques et des côtés des polygones avec leurs longueurs respectives.

\section{Topologie et géométrie des surfaces de caractéristique -1}

Pour saisir la topologie et la géométrie des surfaces non orientables il
est commode de les concevoir comme des surfaces orientables auxquelles
on ajoute des bords que l'on auto-recolle. Dans le cas qui nous intéresse,
nous appuyerons nos raisonnements sur deux représentations~:
une surface$_{-1}$ s'identifie à un tore avec un bord auto-recollé ou
à un pantalon dont les trois bords sont auto-recollés. Ce sont les seules
représentations de ce type.

\subsection{Ovale d'une surface$_{-1}$}
Nous remarquons ici l'existence d'une courbe canonique, à homotopie près, dans
les surface sans bord de caractéristique -1. La géodésique associée
jouera un rôle géométrique essentiel.

\begin{thm-defi}
Chaque surface$_{-1}$ $X$ possède une unique géodésique se relevant
en une géodésique séparante de $X^{or}$. Nous l'appellerons \emph{ovale}
de la surface $X$ et la noterons $\g_X$.
\end{thm-defi}

\begin{figure}[ht]
\centering
\psfrag{T}{$\T_X$}
\psfrag{g}{$\gamma$}
\psfrag{h}{$\delta$}
\psfrag{o}{$\gamma_X$}

\psfrag{X}{$X^{or}$}
\psfrag{g1}{$\gamma^{or}_1$}
\psfrag{g2}{$\gamma^{or}_2$}
\psfrag{h1}{$\delta^{or}_1$}
\psfrag{o1}{$\gamma^{or}_X$}
\psfrag{h2}{$\delta^{or}_2$}
\psfrag{gm}{$\gamma^{or}_2$}

\includegraphics[width=6cm,keepaspectratio=true]{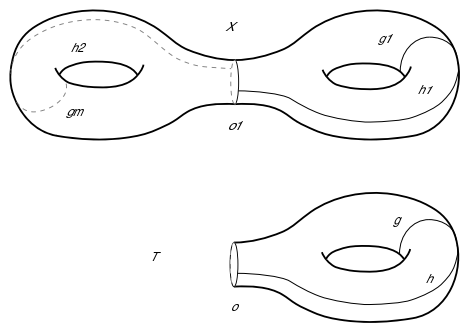}
\caption{Revêtement des orientations}
\label{orientations}
\end{figure}

Dans la première partie de la démonstration de ce résultat, nous exposons une construction élémentaire du revêtement des orientations~: pour une surface$_{-1}$ $X$ nous prenons deux copies
$\T_+$ et $\T_-$ du tore $\T_X$ issu de la découpe de $\g_X$ dans $X$, nous les collons suivant leurs bords en faisant glisser l'un d'un demi-tour par rapport à l'autre (voir figure~\ref{orientations}). Le quotient de la surface obtenue par la translation-reflexion
d'odre 2 le long de la géodésique associée aux bords s'identifie au quotient de $\T_X$ par l'auto-recollement du bord, c'est-à-dire à $X$.
\medskip

\begin{dem} 
Prenons deux tores avec un bord isométriques, nous les collons selon leurs bords en faisant glisser
l'un d'un demi tour par rapport à l'autre. La surface de genre deux ainsi
construite possède une involution anti-holomorphe~: la translation-reflexion d'odre 2 le long du 
bord de chacun des tores. Le quotient de la surface de genre deux par cette involution est une surface$_{-1}$,
en fait la surface de genre deux est le revêtement des orientations de la surface$_{-1}$.
Le bord des tores se projette sur une géodésique non orientable, le relevé de cette géodésique
au revêtement des orientations est une géodésique de longueur double, évidemment séparante.
Nous venons de prouver l'existence, dans une surface$_{-1}$, d'une géodésique (non orientable)
dont le relevé est séparant. 

 Soit $X$ une surface$_{-1}$ et $\g_1$ une géodésique de $X$ se relevant
en une géodésique séparante. Tout d'abord $\g_1$ se relève de manière unique puisque
deux géodésiques séparantes s'intersectent. Ensuite $\g_1$ est non orientable car l'involution
anti-holomorphe de $X^{or}$ la stabilise et agit sur elle sans point fixe.
Montrons son unicité. Supposons l'existence d'une autre géodésique $\g_2$ répondant aux mêmes critères.
Les relevés $\g^{or}_1$ et $\g^{or}_2$ s'intersectent nécessairement.
Découpons $\g^{or}_1$ dans $X^{or}$, nous avons deux tores et $X$ s'obtient par 
l'auto-recollement du bord de l'un d'eux.
La non-orientabilité de $\g_1$ et $\g_2$ implique que le cardinal de $\g_1\cap\g_2$ est impair,
en particulier le nombre de composantes de $\g^{or}_2$ dans chacun des tores
bordés par $\g^{or}_1$ est impair. Ainsi l'involution hyperelliptique de $X^{or}$ va fixer
une de ces composantes qui, nécessairement,
passera par un point de Weierstrass~; d'où $\g^{or}_2$ passe par un point
de Weierstrass. C'est absurde puisqu'elle est séparante.
\end{dem}\hfill$\square$

\begin{rem}
L'ovale est une géodésique non orientable
\end{rem}

\begin{rem}
L'élément du groupe fondamental associé à l'ovale joue un rôle similaire à celui associé au bord
pour le tore à un bord. Un automorphisme du groupe fondamental l'enverra sur lui même, à conjugaison et
inverse près.
\end{rem}

 De par sa caractérisation topologique, l'ovale $\g_X$ d'une surface$_{-1}$ $X$ est la seule géodésique à produire un tore à un bord par sa découpe. Nous noterons $\T_X$ ce tore.
L'application qui à $X$ associe $\T_X$ établit une bijection entre l'espace des surfaces$_{-1}$
et celui des tores à un bord hyperboliques.
\medskip

\subsection{Involution hyperelliptique}\label{partie involution}
Soit $X$ une surface fermée de caractéristique -1, la surface de genre deux $X^{or}$ et le tore $\T_X$ admettent tous deux une involution hyperelliptique. Nous renvoyons à \cite{haas} pour leur étude et observons~:

\begin{pro-defi}
Les involutions hyperelliptiques de $X^{or}$ et $\T_X$ induisent une même involution sur
$X$. Nous l'appellerons \emph{involution hyperelliptique} et la noterons $\iota_X$.

L'ensemble de ses points fixes est formé de l'ovale $\g_X$ et de trois points isolés $w_1$, $w_2$ et $w_3$,
que nous nommerons \emph{points de Weierstrass}.
\end{pro-defi}

\begin{dem}
L'involution hyperelliptique $\iota_{\T_X}$ agit par translation d'ordre 2
sur le bord, ainsi $\iota_{\T_X}$ passe au quotient après auto-recollement.
De même, l'involution $\iota_{X^{or}}$ agit par translation d'ordre 2 sur les géodésiques
séparantes, en particulier elle commute avec la translation-reflexion d'odre 2 le long de $\g^{or}_X$,
et par conséquent passe au quotient.

Une géodésique séparante d'une surface de genre deux borde deux tores, la restriction de l'involution hyperelliptique
à un de ces tores est l'involution hyperelliptique de ce tore. Nous concluons, au vu de la construction du revêtement
des orientations, que les involutions hyperelliptiques de $\T_X$ et $X^{or}$ induisent la même involution de $X$.

L'ensemble des points fixes de $\iota_X$ se déduit trivialement de celui de $\iota_{\T_X}$.
\end{dem}\hfill $\square$\\

 Plusieurs propriétés de l'involution hyperelliptique, bien connues dans le cas orientable, glissent trivialement au cas non orientable via le revêtement des orientations. Soit $X$ une surface fermée de caractéristique -1,
 
\begin{cor}
L'involution hyperelliptique $\iota_X$ est un élément du centre du groupe des automorphismes $\Aut(X)$.
\end{cor}

\begin{cor}
L'involution hyperelliptique stabilise toutes les géodésiques.
\end{cor}

Nous savons depuis les travaux de W.~Scherrer (voir \cite{bujalance2}) qu'une involution d'une surface
non orientable compacte sans bord vérifie l'inégalité $r+2s\leq g$ où $r$ est le nombre de points fixes isolés,
$s$ le nombre de géodésiques de points fixes et $g$ le genre topologique de la surface.
Il en ressort que l'involution hyperelliptique d'une surface$_{-1}$ est l'involution ayant le plus de points fixes isolés.
Dans la partie 2 de notre article, nous étudierons tous les groupes d'automorphismes de surfaces$_{-1}$~;
à cette occasion nous verrons~:

\begin{pro}
Un automorphisme d'une surface$_{-1}$ possède au plus trois points fixes isolés et
l'involution hyperelliptique est le seul à atteindre cette borne.
\end{pro}

Définissons l'involution hyperelliptique d'un pantalon hyperbolique comme la reflexion
par rapport aux perpendiculaires communes de ses bords. Nous avons le corollaire suivant~:

\begin{cor}
L'involution hyperelliptique d'un pantalon hyperbolique passe au quotient après auto-recollement des bords en
l'involution hyperelliptique d'une surface hyperbolique fermée de caractéristique -1.
\end{cor}

Soit $P$ un pantalon hyperbolique. Par auto-recollement d'un, deux ou trois bords de $P$ nous
construisons respectivement un plan projectif à deux bords $Z$, une bouteille de Klein à un bord $Y$
et une surface$_{-1}$ $X$. L'involution hyperelliptique de $P$ passe au quotient après auto-recollement d'un ou plusieurs bords. 
Les involutions qui s'en déduisent seront encore qualifiées d'\emph{hyperelliptiques} et notées $\iota_Y$, $\iota_Z$, $\iota_X$~; d'une manière générale nous préserverons le vocable et les notations définis dans le cadre des surfaces$_{-1}$.

\subsection{Géodésiques}

Nous allons classer très simplement les géodésiques des surfaces$_{-1}$. Mais avant, soulignons
un phénomène assez étonnant.

\begin{lem} 
Dans une surface fermée de caractéristique -1, toute intersection de l'ovale avec une géodésique fermée simple s'effectue à angle droit.
\end{lem}

\begin{dem}
Une géodésique fermée simple est stable par l'involution hyperelliptique, or cette dernière est la reflexion par rapport à l'ovale~; d'où l'angle droit.\hfill$\square$
\end{dem}

\begin{pro}
Il y a trois types de géodésiques fermées simples dans une surface hyperbolique fermée de caractéristique -1~:
\begin{itemize}
\item l'ovale,
\item les géodésiques fermées simples orientables, chacune d'elles passe par exactement deux points de Weierstrass
et ne coupe pas l'ovale~;
\item les géodésiques fermées simples non orientables, chacune d'elles passe par exactement un point de Weierstrass
et intersecte l'ovale en un point.
\end{itemize}
\end{pro}

\begin{dem}
Soient $X$ une surface$_{-1}$ et $\g$ une géodésique distincte de l'ovale.
Commençons par montrer que $\g$ intersecte au plus une fois $\g_X$.

 Supposons que $\g$ coupe $\g_X$ en au moins deux points. Vu dans
le tore $\T_X$, $\g$ possède plusieurs composantes connexes, la parité de leur nombre
décidant de son orientabilité. Considérons une composante connexe $\g_0$ de $\g$,
ou bien $\g_0$ est stable par l'involution, ou bien elle est envoyée sur une autre
composante connexe. Si elle est stable, alors les points d'intersection
de $\g_0$ avec $\partial\T_X$ sont diamétralement opposés et donc identifiés dans $X$~;
donc $\g_0$ est une géodésique simple fermée et $\g=\g_0$, absurde.
Si elle est envoyée sur une autre composante connexe $\g_1$, alors les extrémités
de $\g_1$ sont les points du bord diamétralement opposés aux extrémités de $\g_0$.
Il vient que $\g_0 \cup \g_1$ forme une seule et même géodésique, donc $\g=\g_0\cup\g_1$.
Nous concluons qu'une géodésique de $X$ a au plus deux points d'intersection avec l'ovale.

 Supposons que $\g$ coupe $\g_X$ en exactement deux points, alors $\g$ est orientable. Par ailleurs,
$\g$ ne passe par aucun point de Weierstrass puisque ses composantes $\g_1$ et $\g_2$ sont
échangées par l'involution hyperelliptique. Donc $\g$ se relève en une géodésique séparante sur
laquelle agit l'involution antiholomorphe de $X^{or}$. L'action dans un voisinage de la géodésique est celle d'une
translation-reflexion d'ordre 2 (action sans point fixe), nous en déduisons que $\g$ n'est pas orientable, contradiction~!
D'où $\g$ ne peut couper $\g_X$ en deux points.

 Supposons que $\g$ coupe $\g_X$ en un point, alors $\g$ est non orientable. Comme l'involution hyperelliptique fixe
le point d'intersection elle en fixe nécessairement un autre. Donc $\g$ passe par un point de Weierstrass.

 Supposons que $\g$ n'intersecte pas $\g_X$, alors $\g$ est une géodésique interne de $\T_X$, elle est donc
orientable et passe par deux points de Weierstrass.\hfill$\square$
\end{dem}

\subsection{Dualité}
\begin{figure}[ph]
\centering
\psfrag{T}{$\T_X$}
\psfrag{g}{$\gamma$}
\psfrag{g1}{$\gamma'$}
\psfrag{o}{$\gamma_X$}

\includegraphics[width=6cm,keepaspectratio=true]{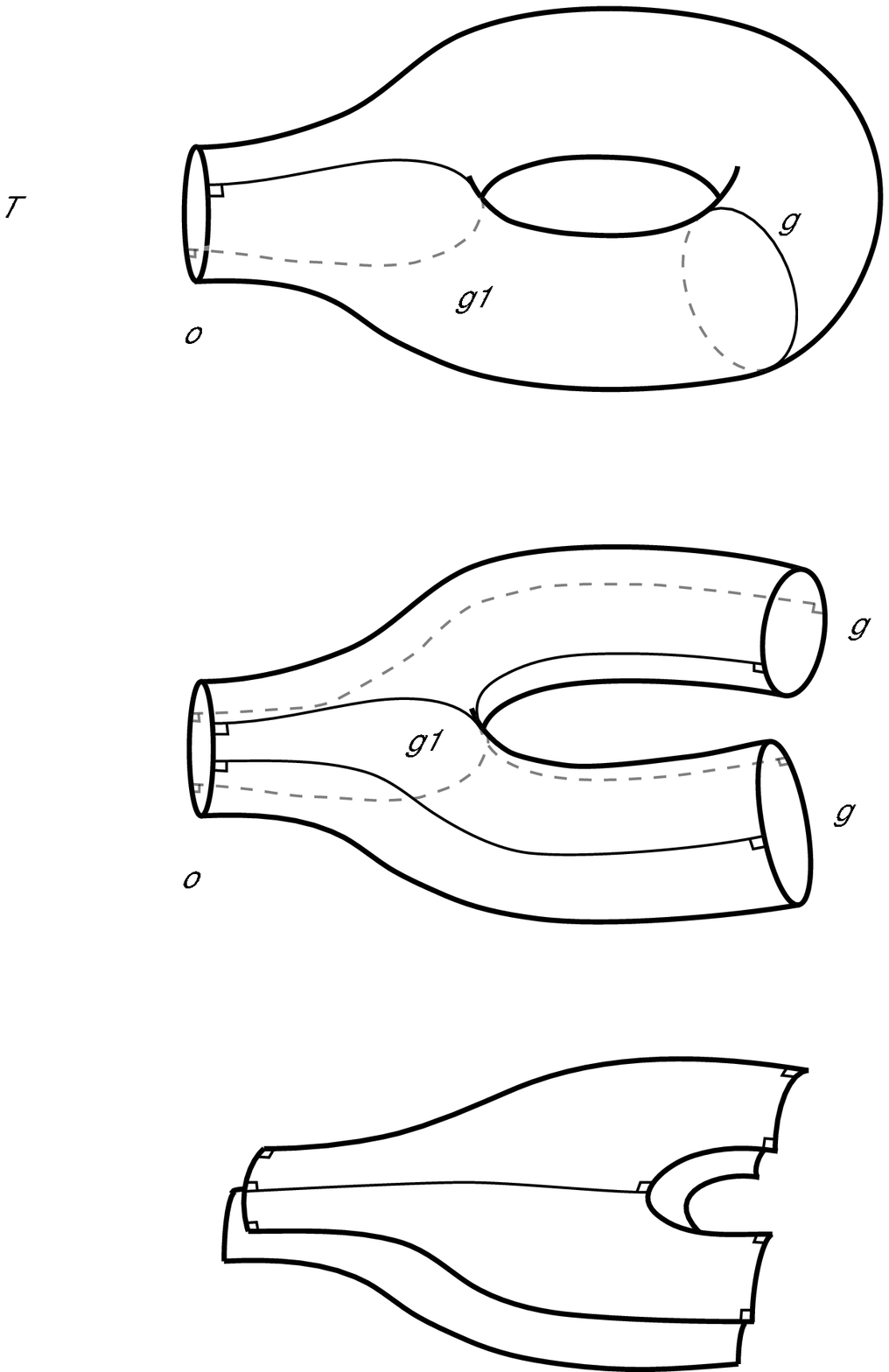}
\caption{Géodésique duale}
\label{dualite}
\vspace{0.5cm}
\includegraphics[width=6cm,keepaspectratio=true]{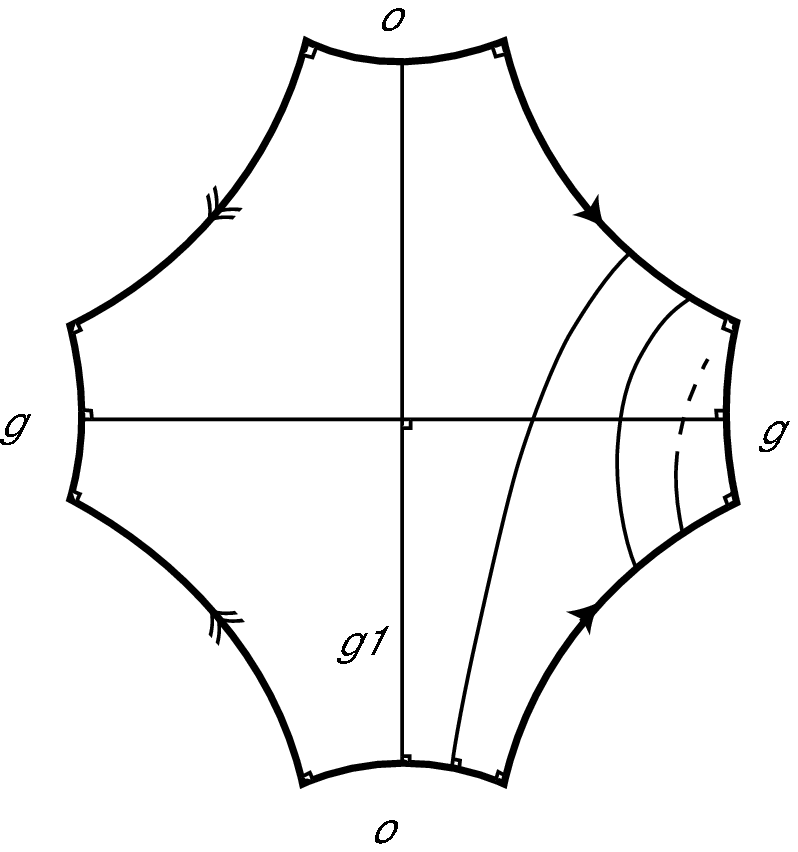}
\caption{}
\label{dualite2}

\end{figure}

Etant donnée une géodésique (fermée simple) orientable $\g$ d'une surface$_{-1}$ $X$, nous pouvons trouver une unique géodésique (fermée simple) non orientable
$\g'$ qui lui est disjointe. Effectivement, nous nous ramenons à un pantalon avec deux bords égaux
en découpant $\g$ dans $\T_X$~; les bords associés à $\g$ admettent une perpendiculaire commune et sa médiatrice
s'identifie à une géodésique non orientable de $X$, disjointe de $\g$ (voir figure~\ref{dualite}).

\begin{pro-defi}
L'application décrite ci-dessus établit une bijection entre les géodésiques fermées simples orientables
et les géodésiques fermées simples non orientables distinctes de l'ovale. Soit $\g$ une géodésique fermée simple, nous noterons $\g'$
la géodésique lui correspondant par cette application et l'appellerons \emph{géodésique duale de $\g$}.
\end{pro-defi}

\begin{dem}
Soit $\g$ une géodésique orientable de $X$, il n'existe pas deux géodésiques non orientables (distinctes de l'ovale)
disjointes de $\g$. Procédons par découpe. Découpons d'abord $\g$ dans $\T_X$, nous obtenons un pantalon avec
deux bords égaux. Découpons ensuite les deux perpendiculaires communes issues de $\partial\T_X$,
nous arrivons à un octogone dont les côtés opposés sont égaux. 
Dans cet octogone essayons de construire une géodésique non orientable (distincte de $\g'$ et $\g_X$)
disjointe de $\g$, cette géodésique partirait à angle droit d'un côté associé à $\g_X$ elle éviterait les
côtés associés à $\g$ et au moins dans un premier temps ceux associés à $\g'$, elle suivrait donc le tracé de la
figure~\ref{dualite2} , c'est impossible car elle ne passerait pas par un point de Weierstrass.

Soit $\g$ une géodésique non orientable de $X$, il existe une et une seule géodésique orientable disjointe de $\g$.
Découpons $\g$ dans $\T_X$, nous obtenons un cylindre, à homotopie près celui-ci ne contient qu'une seule courbe
fermée simple.\hfill$\square$
\end{dem}

\begin{rem}\label{rem dualite}
La dualité est croissante (pour les longueurs). En effet, entre une géodésique orientable $\g$, sa duale $\g'$ et l'ovale $\gamma_X$, nous avons la relation
\begin{eqnarray}
\cosh(l(\gamma)/2)=\sinh(l(\gamma')/2)\sinh(l(\gamma_X)/2).
\end{eqnarray} 
En particulier, la géodésique duale d'une systole
de $\T_X$ est une géodésique non orientable de $X$ distincte de $\g_X$ de plus petite longueur~; nous pouvons montrer (\cite{bavard2}) qu'elle
réalise le rayon d'injectivité de $\gamma_X$. Nous appellerons ces géodésiques \emph{rayons extrémaux}.
\end{rem}

\subsection{Deux exemples de surface$_{-1}$}

\begin{figure}[h]
\centering
\psfrag{o}{$\gamma_X$}

\includegraphics[width=6cm,keepaspectratio=true]{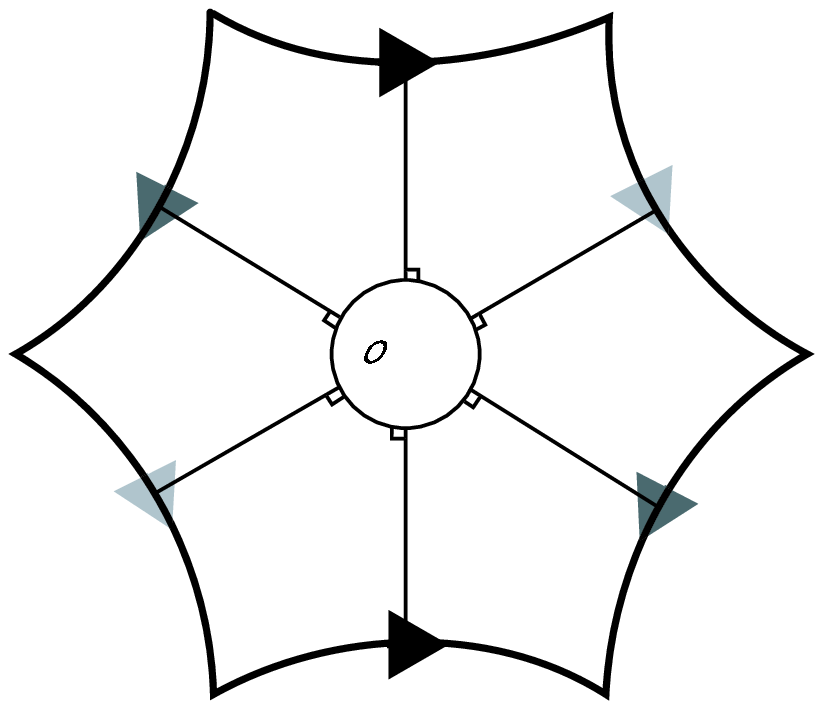}
\caption{Construction du tore équilatéral}
\label{hexagone}

\end{figure}

Nous construisons un tore équilatéral à un bord en identifiant les côtés opposés d'un hexagone régulier, d'angle au sommet
$\pi/3$, auquel on a ôté un disque ouvert (voir figure \ref{hexagone}).
Nous pouvons imposer que la distance du bord aux côtés de l'hexagone soit le quart de la longueur du bord. Nous
noterons $\T_{X(H)}$ le tore ainsi obtenu et $X(H)$ la surface$_{-1}$ s'en déduisant, que nous nommerons \emph{surface
équilatérale}.

 Les points de Weierstrass de $X(H)$ sont les milieux des côtés de l'hexagone. Nous avons représenté les perpendiculaires
communes du bord et des côtés, elles forment trois rayons extrémaux. Les segments joignant les milieux des côtés
constituent les systoles de $\T_X$ puisque ce sont les duales des rayons extrémaux. Il est facile de voir que la longueur des rayons extrémaux est plus petite que
la longueur des systoles de $\T_X$, et nous en déduisons que $X(H)$ possède exactement quatre systoles~: les rayons
extrémaux et l'ovale.\medskip

Une deuxième surface intéressante provient du pentagone droit régulier. Collons deux tels pentagones suivant
un côté nous obtenons un hexagone droit, à partir de deux hexagones nous fabriquons un pantalon, puis en recollant
deux bords de ce pantalon avec un twist nul nous parvenons à une surface$_{-1}$ notée $X(P)$ (nous effectuons la démarche
inverse de celle de la figure~\ref{dualite}).

Cette surface possède cinq systole~: l'ovale, deux systoles orientables et deux rayons extrémaux. En fait, il y
a une systole par côté du pentagone. Il est très facile de voir que ces géodésiques sont bien les systoles, car toute géodésique passe par un sommet d'un pentagone (un point de Weierstrass) avant de traverser celui-ci.

\subsection{La bouteille de Klein à un bord et le plan projectif à deux bords}
Nous reprenons les notations de la fin du paragraphe \ref{partie involution}.
Le cas des surfaces$_{-1}$ a visiblement valeur de paradigme.
En raisonnant avec $Y$ et $Z$ comme nous le ferions avec $X$ nous trouvons~:

\begin{pro}
La bouteille de Klein à un bord admet une seule géodésique fermée simple orientable et une infinité de géodésiques fermées simples non orientables.
\end{pro}

\begin{dem}
La géodésique orientable est la duale du bord. Des exemples de géodésiques non orientables sont les 
deux bords auto-recollés de $P$, nous obtenons une infinité d'autres géodésiques non orientables en faisant agir sur elles
le twist selon la géodésique orientable.\hfill$\square$
\end{dem}

\begin{rem}
En fait, nous récupérons toutes les géodésiques non orientables par la méthode explicitée dans la preuve.
\end{rem}

\begin{pro}\label{geodesiques_plan_projectif}
Le plan projectif à deux bords n'admet aucune géodésique fermée simple orientable et seulement deux géodésiques fermée simples non orientables.
\end{pro}

\begin{dem}
Une géodésique orientable serait duale de chacun des bords, c'est impossible. 
Il y a une géodésique non orientable évidente~: le bord auto-recollé de $P$.
Cherchons les autres géodésiques non orientables. Elles sont orthogonales à l'ovale et l'intersecte
en un seul point, de plus elles passe obligatoirement par les point de Weierstrass de $Z$.
Nous voyons rapidement qu'il n'y a q'une possibilité.\hfill$\square$
\end{dem}

\begin{figure}[ht]
\centering
\psfrag{g}{$\gamma$}
\psfrag{g1}{$\gamma_1$}
\psfrag{g2}{$\gamma_2$}
\psfrag{o}{$\gamma_X$}
\psfrag{tg}{$t_{\gamma}(\gamma_1)$}
\psfrag{w1}{$w_1$}
\psfrag{w2}{$w_2$}

\includegraphics[width=6cm,keepaspectratio=true]{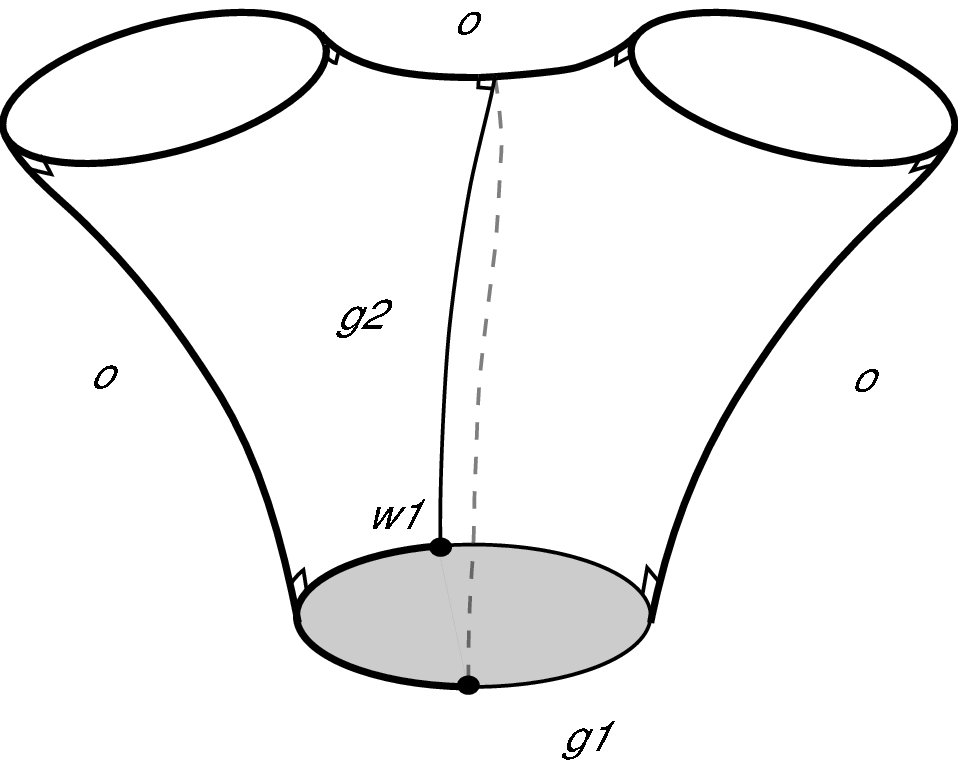}
\caption{Le plan projectif à deux bords}
\label{plan projectif}
\vspace{0.5cm}
\includegraphics[width=6cm,keepaspectratio=true]{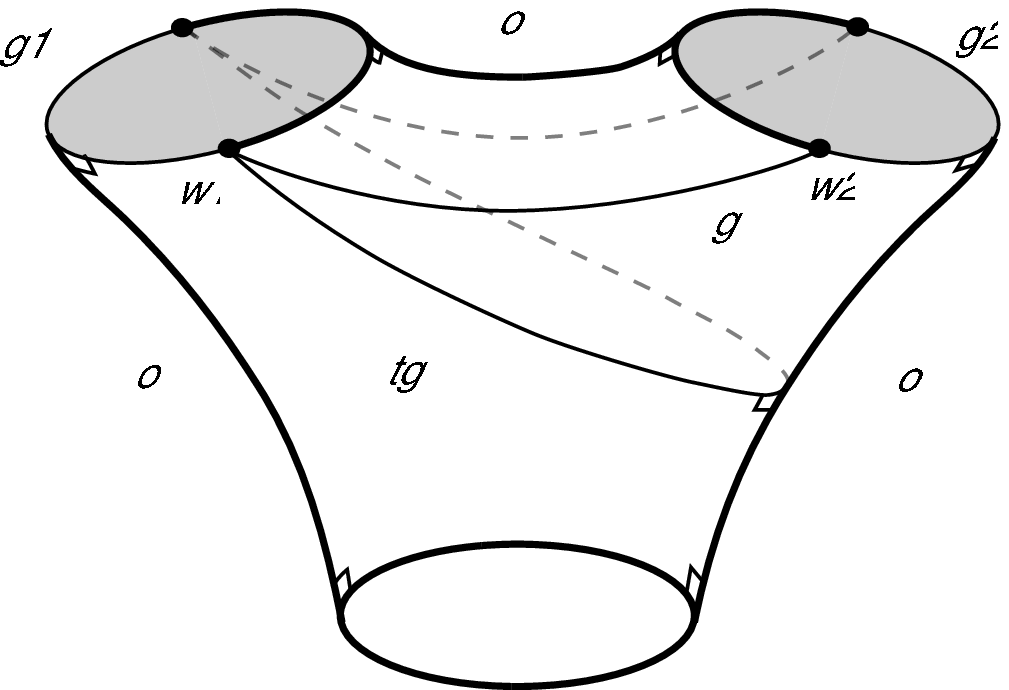}
\caption{La bouteille de Klein à bord}
\label{bouteille}

\end{figure}

 Pour chacune des figures ci-dessus, les bords grisés du pantalon représentent ceux auto-recollés.
Sur la figure~\ref{plan projectif} nous avons noté $\g_1$ et $\g_2$
les deux géodésiques du plan projectif à deux bords. Sur la figure~\ref{bouteille} nous avons
noté $\gamma$ la géodésique orientable et $\gamma_1$, $\gamma_2$ deux géodésiques non orientables~; le $t_\g$ désigne le twist selon $\g$.

\subsection{Commentaires}
Dans notre approche des surfaces non orientables nous avons adopté un point de vue métrique, hyperbolique. Un point de vue plus classique serait celui
des surfaces de Klein, c'est-à-dire des surfaces munies d'un atlas dont les changements de cartes sont holomorphes
ou anti-holomorphes. L'étude des surfaces de Klein a connu un fort développement depuis l'impulsion
donnée par le livre \cite{alling} d'Alling et Greenleaf au début des années 70~;
les méthode utilisées, basées sur l'emploie des groupes non euclidiens cristallographique
(\guillemotleft N.E.C groups \guillemotright\ en anglais)
font essentiellement appel à des techniques combinatoires. Les résultats obtenus concernent
principalement les groupes d'automorphismes (borne du nombre d'automorphismes, détermination des groupes
d'automorphismes des surfaces hyperelliptiques...). Dans ce cadre l'involution hyperelliptique est
souvent définie comme une involution $\s$ d'une surface $X$ telle que $X/\langle \s \rangle$ soit de genre algébrique
0, il est alors montré que l'involution hyperelliptique est unique et dans le centre de $\Aut(X)$ (voir l'article \cite{jbujalance} de J.A. Bujalance).

 Parmi les différents articles de E. Bujalance, J.J. Etayo, J.M. Gamboa, G. Gromadzki,
et D. Singerman, nous renvoyons à \cite{bujalance3}
pour une introduction aux surfaces de Klein, à \cite{bujalance4} et \cite{jbujalance} pour l'étude de l'involution hyperelliptique des surfaces de Klein sans bord. 

\section{Actions de groupes d'automorphismes}

Commençons par rappeler un résultat connu sur le sujet~; c'est un cas particulier du
théorème principal de \cite{bujalance1}.

\begin{pro} Les seuls groupes se réalisant comme groupe des automorphismes d'une surface$_{-1}$
sont~: $\Z/2\Z$, $\Z/2\Z \times \Z/2\Z$, $\D_4$ et $\D_6.$
\end{pro}

 Nous allons décrire toutes les surfaces$_{-1}$ possédant un groupe d'automorphismes différent de $\Z/2\Z$, mais
cela passe d'abord par l'étude des involutions.\medskip

 Dans toute la suite $X$ désignera une surface$_{-1}$. Soulignons que le paramètre de twist est à valeurs dans le cercle $\R/\Z$.

\subsection{Actions de $\Z/2\Z\times \Z/2\Z$}
Dans \cite{bujalance2}, nous trouvons des informations précises, attribuées à Scherrer,
sur l'action par involution. Pour une involution d'une surface$_{-1}$, notons $r$ son nombre de points fixes
isolés, $s^+$ (resp. $s^-$) le nombre de géodésiques orientables (resp. non orientables) fixées points à points.
Alors

\begin{pro} Les seuls triplets $(r,s^+,s^-)$ réalisés par l'involution d'une surface$_{-1}$
sont~: $(1,0,1)$, $(1,1,1)$ et $(3,0,1)$. Et à chacun d'eux correspond une unique involution, à conjugaison
topologique près.
\end{pro}

Procédons à quelques remarques. Soit $\s$ une involution d'une surface$_{-1}$ $X$,
nous noterons $\d^-$ sa géodésique non orientable de points fixes et $\d^+$
la géodésique duale de cette dernière.
\begin{itemize}
\item Si $(r,s^+,s^-)=(3,0,1)$ alors $\s=\iota_X$.
En effet, nous avons vu que $\s$ est conjuguée à $\iota_X$,
mais la classe d'homotopie de $\g_X$ est fixe par homéomorphisme, d'où $\d^-=\g_X$ et $\s=\iota_X$.
\item Si $\s\neq \iota_X$, alors $\s$ a exactement deux points fixes sur $\g_X$ dont un appartenant à $\d^-$.
\item Si $\s\neq \iota_X$, alors $\s$ peut agir de deux manières sur $\d^+$~: ou bien $\s$ fixe tous les points
et $(r,s^+,s^-)=(1,1,1)$, 
ou bien $\s$ agit par translation-reflexion d'ordre 2 et $(r,s^+,s^-)=(1,0,1)$.
En tout cas, $\s$ ne peut agir par translation pour des raisons de connexité évidentes lorsqu'on se place dans $\T_X$.
\end{itemize}

 Considérons maintenant une surface $X$ sur laquelle agit $G=\lbrace id_X, \s_1,\s_2,\iota_X \rbrace$
un groupe d'automorphismes isomorphe à $\Z/2\Z\times \Z/2\Z$. Dans la suite nous adopterons les notations suivantes~:
$\d_i^-$ désignera la géodésique non orientable fixée point à point par $\s_i$ ($i=1,2$), et $\d_i^+$ sa 
géodésique duale.\medskip

Comme $\s_1 \s_2=\iota_X$, il vient facilement $(r,s^+,s^-)_{\s_1}=(r,s^+,s^-)_{\s_2}$, nous noterons désormais
$(r,s^+,s^-)_G$ ce triplet. Nous allons inventorier les actions possibles de $G$ selon les deux valeurs éventuelles
de ce triplet.\medskip

\subsubsection{Si $(r,s^+,s^-)_G=(1,1,1)$} Alors $\d_1^-$ et $\d_2^-$ sont disjointes, sinon elles se
couperaient en un point de Weierstrass et, $\d_1^-$ et $\d_2^+$ s'intersecteraient en dehors des points
de Weierstrass ce qui est impossible puisque ce sont les points fixes de $\iota_X=\s_1 \s_2$.
Nous en déduisons
$$\d_1^+\cap\d_2^+=\lbrace w_1 \rbrace,\quad \d_1^+\cap\d_2^-=\lbrace w_2\rbrace,\quad \d_1^-\cap\d_2^+=\lbrace w_3 \rbrace,\quad \d_1^-\cap\d_2^-=\emptyset, $$
où $w_1,\ w_2$ et $w_3$ sont les points de Weierstrass de $X$. Souligons que chaque intersection non vide
s'accomplie avec un angle droit.\medskip

\begin{figure}[ht]

\psfrag{T}{$\T_X$}
\psfrag{o}{$\gamma_X$}
\psfrag{w2}{$w_2$}
\psfrag{w1}{$w_1$}
\psfrag{w3}{$w_3$}

\psfrag{d1}{$\delta^+_1$}
\psfrag{n1}{$\delta^-_1$}
\psfrag{d2}{$\delta^+_2$}
\psfrag{n2}{$\delta^-_2$}

\begin{minipage}[c]{0.5\linewidth}
\centering
\includegraphics[width=6cm,keepaspectratio=true]{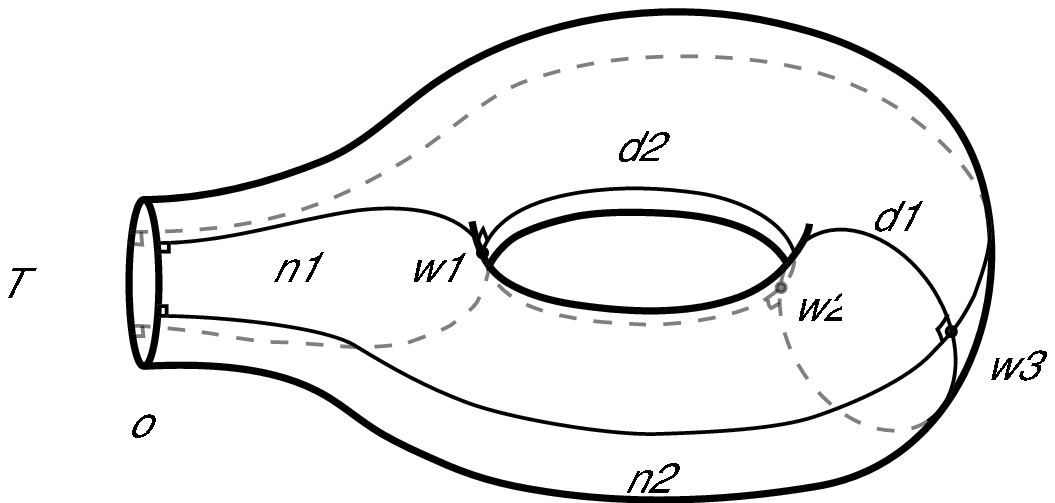}
\caption{Surface avec un twist nul selon $\delta_i^+$}
\label{twist_0}
\end{minipage}%
\vspace{0.5cm}
\begin{minipage}[c]{0.5\linewidth}
\centering
\includegraphics[width=6cm,keepaspectratio=true]{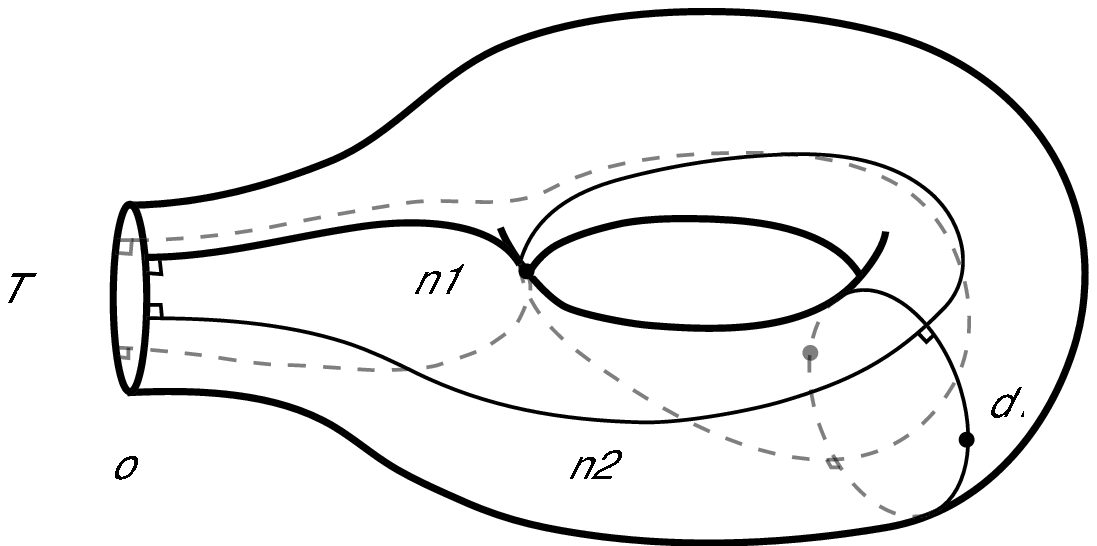}
\caption{Surface avec un twist 1/2 selon $\delta_i^+$}
\label{twist_demi}
\end{minipage}%

\end{figure}

 Contemplons $X$, nous voyons que les coefficients de twist relatifs aux $\d_i^+$ ($i=1,2$) sont nuls
(voir figure~\ref{twist_0}). Réciproquement,
si $X$ est une surface dont le twist selon une certaine géodésique intérieure est nul, alors il est facile de construire
un groupe d'automorphismes $G\simeq \Z/2\Z \times \Z/2\Z$ contenant l'involution hyperelliptique
tel que $(r,s^+,s^-)_G=(1,1,1).$\medskip

\subsubsection{Si $(r,s^+,s^-)_G=(1,0,1)$} Montrons que $\d_1^+\cap \d_2^+=\lbrace w_2,w_3 \rbrace$.
Premièrement $\d_1^+$ et $\d_2^+$ passent par les mêmes points de Weierstrass~; dans le cas contraire $\d_1^-$ et $\d_2^+$
s'intersecteraient en un point de Weierstrass, et ce dernier ne serait pas fixé par $\iota_X$, c'est impossible.
Ensuite $\d_2^+$ et $\d_1^+$ sont toutes deux stables par $\s_1$ et $\s_2$ qui agissent
par translation-reflexion sur elles. L'existence d'un
troisième point d'intersection contredirait l'égalité $\iota_X=\s_1\s_2$, il suffit de regarder comment
$\s_1$, $\s_2$ et $\iota_X$ déplaceraient ce point sur $\d_1^+$. Donc
 $$\d_1^-\cap\d_2^-=\lbrace w_1 \rbrace\quad \rm{et}\quad \d_1^+\cap \d_2^+=\lbrace w_2,w_3 \rbrace.$$ 
En découpant $\d_1^-$ et $\d_1^+$ dans $\T_X$, nous devinons que 
$\d_1^-$ et $\d_2^+$ s'intersectent en exactement deux points $\lbrace p_1,p_2 \rbrace.$
Bien sûr, $p_1$ et $p_2$ sont échangés par $\s_2$, et l'intersection en ces points comme celle en $w_1$
s'effectue avec un angle droit car $\s_1$ (la reflexion suivant $\d_1^-$) stabilise $\d_2^+$.
Nous concluons que le paramètre de twist
de $X$ selon $\d_i^+$ ($i=1$ ou $2$) vaut 1/2 (voir figure~\ref{twist_demi}).
Inversement si $X$ possède un paramètre de twist de 1/2 selon
une certaine géodésique interne de $\T_X$, alors nous pouvons formé un groupe d'automorphisme
$G\simeq \Z/2\Z \times \Z/2\Z$ contenant l'involution hyperelliptique et vérifiant $(r,s^+,s^-)_G=(1,0,1).$\medskip

 Les familles de surfaces dont nous venons de parler constituent des lieus de points fixes relatifs à certains
sous-groupes du groupe modulaire. Vu dans l'espace de Teichmüller elles formeront des sous-variétés totalement
géodésiques. Ceci sera expliciter plus loin.

\subsection{Action de $\D_4$}
 Dans ce paragraphe nous supposons $\Aut(X)\simeq \D_4$. Soit $\tau$ un
automorphisme d'ordre 4 de $X$, $\tau$ agit par translation d'odre 4 sur $\partial\T_X$ (quitte à le multiplier
par $\iota_X$).\medskip

\begin{figure}[h]

\psfrag{a}{$a$}
\psfrag{b}{$b$}
\psfrag{c}{$c$}
\psfrag{d}{$d$}
\psfrag{e}{$e$}

\centering
\includegraphics[width=6cm,keepaspectratio=true]{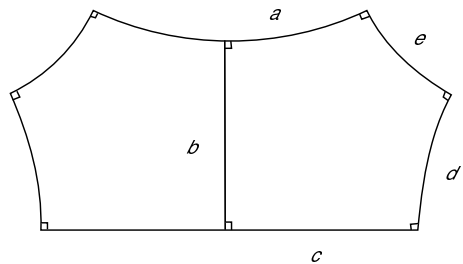}
\caption{Hexagone et Pentagone droits}
\label{pentagone1}

\end{figure}

 Soit $s_1$ une systole de $\T_X$, $s_1$ passe par deux points de Weierstrass.
Clairement $\tau$ fixe un point de Weierstrass, disons $w_1$, et échange les deux autres.
Supposons que $s_1$ ne passe pas par $w_1$, obligatoirement $\tau$ stabilisera $s_1$ car deux systoles ne
peuvent s'intersecter plus d'une fois. Cependant ceci impliquerait que $\tau^2=\iota_X$ fixe point à point
$s_1$, absurde. Donc $s_1$ passe par $w_1$ et par un deuxième point de Weierstrass, par exemple $w_2$.
En conséquence, $\tau$ étant la rotation d'angle $\pi/2$ en $w_1$ envoie $s_1$ sur une deuxième systole
$s_2$ qui lui est orthogonale. En découpant $s_1$ dans $\T_X$, nous obtenons un pantalon dont la perpendiculaire
commune aux bords identifiés à $s_1$ est $s_2$~; le twist de $X$ par rapport à $s_1$ est de ce fait nul.
En outre, en se plaçant dans un pentagone et avec les notations de la figure~\ref{pentagone1}
nous avons $a=e$ et $b=d$.\medskip

 Inversement, partant d'un pentagone vérifiant $a=e$ et $b=d$, et en recollant le pantalon associé avec un twist nul,
nous aboutissons à une surface ayant $\D_4$ comme groupe des automorphismes. Via un argument de trigonométrie
hyperbolique, nous en concluons qu'il y a exactement une famille à un paramètre de surfaces ayant $\D_4$ pour groupes
des automorphismes, et le paramètre est $l(\g_X)$.\medskip

 Quelques constatations sur ces surfaces. Pour chacune de ces surfaces, $\T_X$ a deux systoles et l'ovale
$\g_X$ deux rayons extrémaux.
Trivialement, une seule de ces surfaces possède plus de deux systoles~: la surface construite
à partir du pentagone régulier. Enfin, nous voyons sans peine comment agisse tous les automorphismes~; en particulier
pour un telle surface $X$, $\Aut(X)$ contient deux sous-groupes $G_1$, $G_2$ de la forme
$\lbrace id_X, \s_1,\s_2,\iota_X \rbrace \simeq \Z/2\Z\times\Z/2\Z$ avec $(r,s^+,s^-)_{G_1}=(1,0,1)$ et
$(r,s^+,s^-)_{G_2}=(1,1,1)$ (pour ce dernier nous avons $\d_i^+=s_i$). Cette famille de surfaces
se situe donc à l'intersection des deux familles étudiées précédemment. En fait, elles forment exactement
l'intersection car, nous allons le voir, les surfaces $X$ avec $\Aut(X)\simeq \D_6$ n'ont pas de groupe d'automorphisme
$G\simeq\Z/2\Z \times \Z/2\Z$, contenant $\iota_X$, avec $(r,s^+,s^-)_{G}=(1,1,1)$.

\subsection{Action de $\D_6$} Soit $X$ une surface avec $\Aut(X)\simeq \D_6$.
Considérons $\tau$ un automorphisme d'ordre 6 de $X$, nous avons $\tau^3=\iota_X$. Vu que $\iota_X$ fixe les points
de Weierstrass, la restriction de $\tau$ à ces points est un cycle d'odre 3 (ce ne peut être l'identité, sinon nous
pourrions par exemple construire des systoles s'intersectant en deux points).\medskip

 Prenons un rayon extrémal de $X$, $\tau$ transforme ce rayon extrémal en un autre rayon extrémal puisqu'il
agit transitivement sur les points de Weierstrass, il en découle que $X$ possède trois rayons extrémaux disjoints
$r_1$, $r_2=\tau(r_1)$ et $r_3=\tau^2(r_1)$.
Découpons ces trois rayons dans $X$, nous obtenons un pantalon avec trois bords de mêmes longueurs.
Les systoles internes de $\T_X$ sont les géodésiques duales des rayons extrémaux, on les visualise avec facilités
sur le pantalon~: elles joignent les milieux de deux bords. Toujours dans le même pantalon, regardons la perpendiculaire
à $\g_X$ issue de $w_1$. Cette perpendiculaire forme une
géodésique fermée simple, orthogonale à la systole $s_1$ et au rayon $r_1$. Si nous regardions cette
perpendiculaire dans le pantalon produit par la découpe de $\g_X$ et $s_1$, alors nous verrions que le twist de $X$
selon $s_1$ a pour valeur 1/2.\medskip

 Réciproquement, en partant d'un pantalon avec trois bords de même longueur, et en auto-recollant ces bords
nous obtenons une surface $X$ avec $\Aut(X)\simeq \D_6$. Conclusion~: il y une famille à un paramètre (la longueur
de l'ovale $\g_X$) de surfaces ayant $\D_6$ pour groupe d'automorphismes. Une représentation agréable de ces
surfaces est celle de la figure~\ref{hexagone}. Le groupe d'automorphisme nous apparaît alors avec limpidité, nous
remarquons en particulier que les trois sous-groupes contenant $\iota_X$ et isomorphes au groupe de Klein
vérifient tous $(r,s^+,s^-)_{G}=(1,0,1)$.

\section{Espaces de Teichmüller et systole des surfaces$_{-1}$}

Dans cette partie, toute surface sera supposée fermée. Pour une géodésique orientable $\g_1$,
le twist $\theta_1$ sera croissant si l'on a un glissement d'un bord dans le sens donné par la flèche de la
figure~\ref{pentagone}~; si la figure~\ref{pentagone} représente un twist nul selon $\g_1$, alors la
figure~\ref{lambda} représente un twist compris entre 0 et 1/2. Nous noterons $t_{\g_1}$ le
twist de Dehn selon $\g_1$.

\subsection{Définition}
Fixons nous une surface hyperbolique orientable $R$ de genre $g\geq 2$.
On appelle \emph{surface marquée de genre $g$} un couple
$(S,\f)$ formé d'une surface hyperbolique orientable $S$ de genre $g$ et d'un homéomorphisme préservant l'orientation
$\f:R\rightarrow S$. Sur l'ensemble des surfaces marquées on définit le relation d'équivalence suivante~:
$$(S,\f)\sim  (S',\f')\ \rm{si}\ \f'\circ \f^{-1}\ \rm{est\ isotope\ \grave{a}\ une\ isom\acute{e}trie}.$$
Nous regarderons l'\emph{espace de Teichmüller} $\t_g$ des surfaces hyperboliques orientables 
de genre $g$ comme l'espace des classes
d'équivalences $[S,\f]$ de surfaces marquées. Il est muni d'une structure de variété analytique complexe
(de dimension $3g-3$) ainsi que d'une métrique kählérienne de courbure négative (dite de Weil-Petersson). \medskip

  Sur $\t_g$, agit naturellement le groupe des homéotopies. Son quotient par le
sous-groupe des éléments agissant trivialement est le groupe modulaire,
que nous noterons $Mod_g$.
L'action de $Mod_g$ sur $\t_g$ s'effectue par isométrie pour la métrique de Weil-Petersson.\medskip

 En supprimant la condition relative à l'orientabilité, on étend au cadre non orientable
les notions de surfaces marquées et d'espace de Teichmüller. L'espace de Teichmüller des surfaces
hyperboliques non orientables de genre $g$ sera noté $\t_g^-$, et le groupe modulaire associé $Mod_g^-$.
En considérant le revêtement des orientations plutôt qu'une surface elle-même, l'espace $\t_g^-$ apparaît
comme lieu des points fixes d'une isométrie négative de $\t_{g+1}$, et en tant que tel c'est une sous-variété
totalement géodésique de $\t_{g+1}$ sur laquelle agit $Mod_g^-$ par isométries.\medskip

 Nous allons nous intéresser à $\t_3^-$ l'espace de Teichmüller des surfaces hyperboliques non orientables
de genre $3$. Vu dans $\t_2$ c'est le lieu des points fixes de l'homéotopie se réalisant comme
translation-reflexion d'ordre 2 le long de la géodésique séparante au-dessus de l'ovale (pour les revêtements
de surfaces non orientables).

\subsection{Paramétrisations de $\t_3^-$} 
Nous allons exhiber deux paramétrisations globales de $\t_3^-$.

\begin{pro}
Soit $\g_1$ une géodésique orientable des surfaces$_{-1}$. Soient $(\theta_1, l_1)$ les coordonnées twist-longueur
selon $\g_1$, l'application 
$$\begin{array}{cccll}
\Phi & : & \t_3^- & \longrightarrow & \R\times\R_+^\ast\times\R_+^\ast \\
     &   & X      & \longmapsto     & (\theta_1(X),l_1(X),l(\g_X))     \\
\end{array}$$
est un difféomorphisme.
\end{pro}

\begin{dem}
 Considérons les coordonnées de Fenchel-Nielsen de $\t_2$ associées au système de géodésiques
$\lbrace \tilde{\g}_X,\tilde{\g}_{11}, \tilde{\g}_{12} \rbrace$, où les $\tilde{\g}_{1i}$ ($i=1,2$)
sont les deux relevés de $\g_1$. Il est connu (voir par exemple \cite{imayoshi}) que l'application
$$\begin{array}{cccll}
\phi & : & \t_2 & \longrightarrow & (\R_+^\ast)^3 \times (\R)^3 \\
     &   & X    & \longmapsto     & (l(\tilde{\g}_X),l_{11}(X),l_{12}(X),\theta_{\tilde{\g}_X}(X),
 \theta_{11}(X), \theta_{12}(X)) \\
\end{array}$$
est un système de coordonnées globales sur $\t_2$.

 L'espace $\t_3^-$ s'identifie clairement à l'ensemble des points de $\t_2$ déterminés par les conditions
$$ \left\{
\begin{array}{rll}
\theta_{\tilde{\g}_X}(X) & = & \pi \\
     l_{11}(X) & = & l_{12}(X) \\
\theta_{11}(X) & = & \theta_{12}(X) 
\end{array} \right. .$$ 
Il s'en suit que $(\theta_1,l_1,l(\g_X))$ forme un système de coordonnées globales sur $\t_3^-$.\nolinebreak$\square$
\end{dem}

\begin{rem}
Dans ce système de coordonnées, il semble naturel d'imaginer $\t_3^-$ comme espace fibré au-dessus de $\R^\ast_+$
(la longueur du bord), chaque fibre étant l'espace de Teichmüller du tore avec un bord de longueur fixée.
\end{rem}

 Le deuxième système de coordonnées que nous allons présenter consiste en trois fonctions longueurs.
Nous ne nous appuyons pas sur le revêtement des orientations.

\begin{thm} Soient $\g_4$, $\gamma_5$ et $\gamma_6$ trois géodésiques non orientables
et disjointes. Alors, l'application ayant pour coordonnées les fonctions longueurs associées
$$\begin{array}{cccll}
\Psi & : & \t_3^- & \longrightarrow & (\R_+^\ast)^3  \\
     &   & X      & \longmapsto     & (l_4(X),l_5(X),l_6(X)) \\
\end{array}$$
établit un difféomorphisme entre $\t_3^-$ et $(\R_+^\ast)^3$.
\end{thm}

\begin{figure}[ht]

\psfrag{o}{$\gamma_X$}
\psfrag{w2}{$w_2$}
\psfrag{w1}{$w_1$}
\psfrag{w3}{$w_3$}

\psfrag{g2}{$\gamma_2$}
\psfrag{g1}{$\g_1$}
\psfrag{g4}{$\g_4$}
\psfrag{g5}{$\g_5$}
\psfrag{l}{$\lambda$}
\psfrag{t}{$\theta_1l_1/2$}

\centering
\includegraphics[width=6cm,keepaspectratio=true]{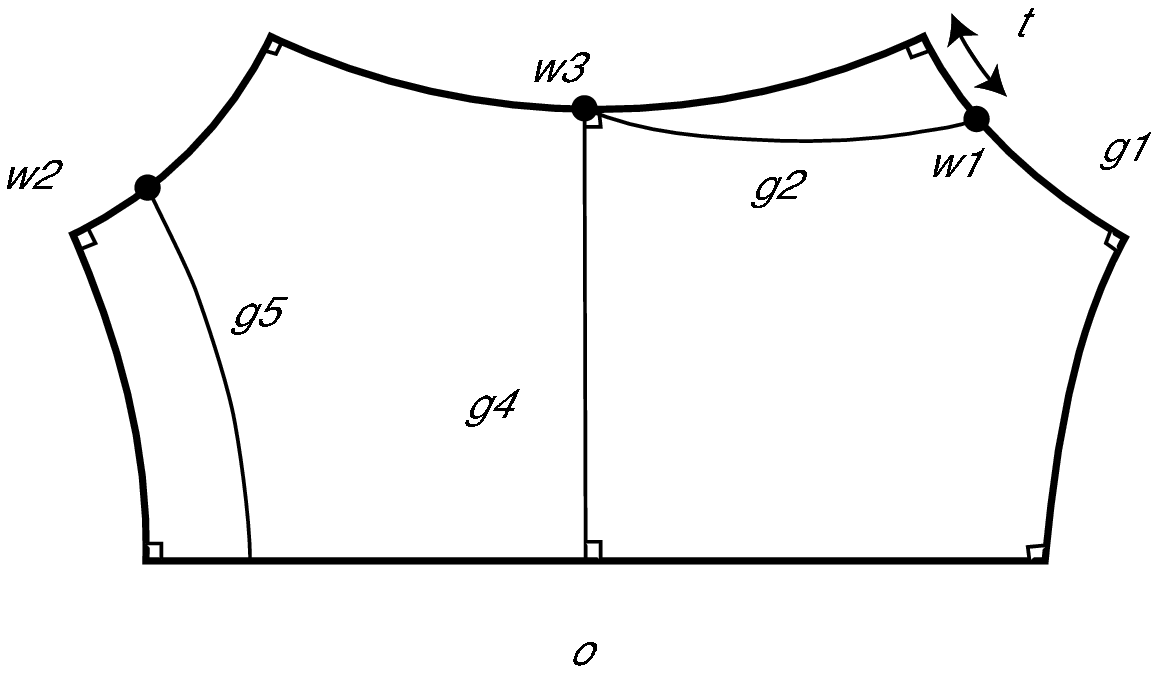}
\caption{}
\label{lambda}
\end{figure}

\begin{dem} Tout d'abord, l'application $\Psi$ est bijective. En effet, l'opération de découpe
des trois géodésiques $\g_4$, $\g_5$ et $\g_6$ établit une bijection entre $\t_3^-$ et l'espace des pantalons~; mais
ce dernier est paramétré par les longueurs des bords, d'où l'affirmation.

 Ensuite, $\Psi$ à l'instar des fonctions longueurs est de classe $C^\infty$ sur $\t_3^-$.
Pour conclure, montrons que
son application réciproque l'est aussi~; autrement dit, montrons que $l(\g_X)$, $l_1$ et $\theta_1$
dépendent de façon $C^\infty$ de $l_4(X)$, $l_5(X)$ et $l_6(X)$. Vu dans le pantalon, $l(\g_X)$ s'exprime comme somme
des longueurs des perpendiculaires communes aux bords, elle est clairement de classe $C^\infty$.
Sans perte de généralités, supposons $\g_4$ et $\g_1$
duales l'une de l'autre. En découpant $\g_1$ dans $\T_X$ (voir figure~\ref{lambda}), nous voyons $l_1$
comme fonction $C^\infty$ en $l_4(X)$ et $l(\g_X)$. Quant-au twist~: introduisons la géodésique $\g_2$
duale de $\g_5$ dans $X$,
à l'aide de relations trigonométriques nous écrivons les égalités
$$\cosh(\theta_1(X)l_1(X)/2)=\cosh(l_2(X)/2)\tanh(l_1(X)/2)\tanh(l_4(X)/2)$$
et (voir Remarque~\ref{rem dualite})
$$\cosh(l_2/2)=\sinh(l_5/2)\sinh(l(\g_X)/2),$$
prouvant ainsi le caractère lisse de $\theta_1$ par rapport à $(l_4(X),l_5(X),l_6(X)).$
\end{dem} $\hfill\square$

\begin{rem} Une paramétrisation de l'espace de Teichmüller par des longueurs est un phènomène très rare,
il ne se produit jamais dans le cas hyperbolique orientable, même locallement (voir par exemple \cite{seppala}).
\end{rem}
\begin{rem} L'espace $\t_3^-$ et l'espace de Teichmüller des pantalons sont difféomorphes.
\end{rem}

\begin{figure}[ht]

\psfrag{o}{$\gamma_X$}
\psfrag{w2}{$w_2$}
\psfrag{w1}{$w_1$}
\psfrag{w3}{$w_3$}

\psfrag{g1}{$\g_1$}
\psfrag{g4}{$\g_4$}
\psfrag{g2}{$\g_2$}
\psfrag{g3}{$\g_3$}

\psfrag{d1}{${\g'}_1$}
\psfrag{d3}{${\g'}_3$}
\psfrag{d2}{${\g'}_2$}

\psfrag{t}{$t$}

\centering
\includegraphics[width=6cm,keepaspectratio=true]{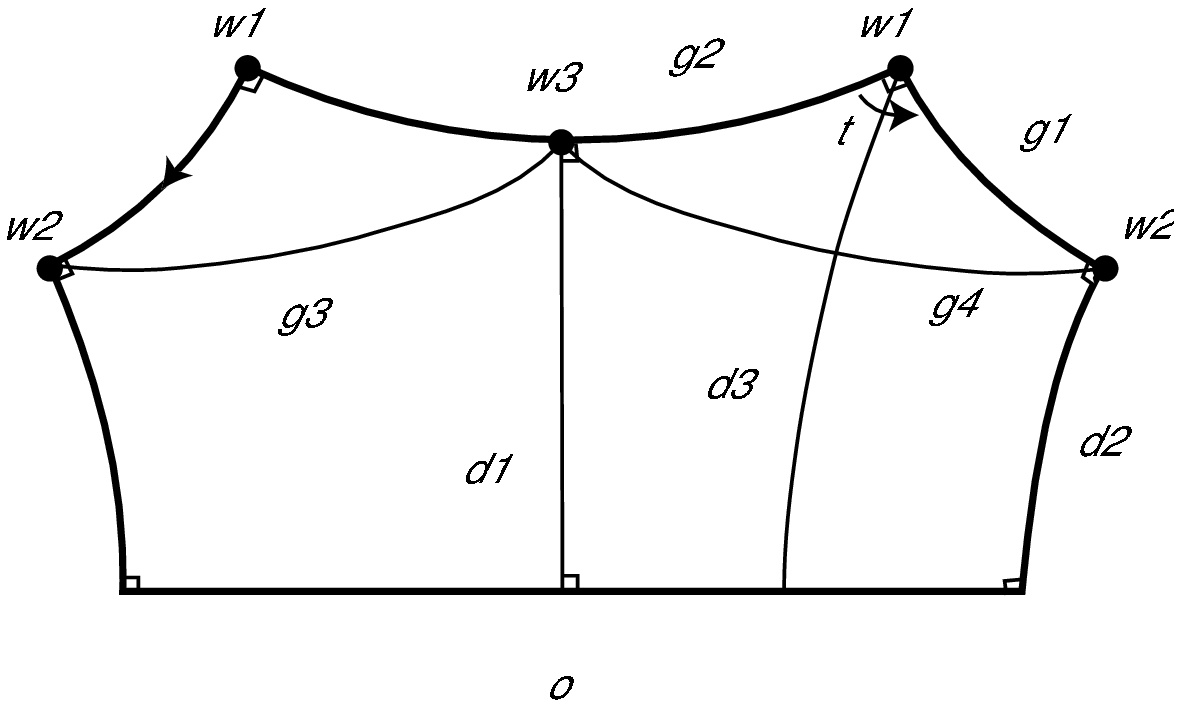}
\caption{}
\label{pentagone}
\end{figure}

 Nous allons nous fixer un système de coordonnées auquel nous nous référerons par la suite.
Choisissons $X(P)$ comme origine du Teichmüller, et adoptons les notations suivantes (voir figure~\ref{pentagone})~:
\begin{itemize}
\item $\g_1$ et $\g_2$ désignent les deux systoles orientables de $X(P)$,
\item $\g_3=t^{-1}_{\g_1}(\g_2)$ et $\g_4=t_{\g_1}(\g_2)$ sont les géodésiques obtenues à partir de $\g_2$
par un seul twist de Dehn selon $\g_1$.
\end{itemize}
Nous optons pour le système de coordonnées $(\theta_1,l_1,l(\g_X))$ avec la convention $\theta_1([X(P),id])=0$.
Nous noterons $[X(H),\f_{1/2}]$ la surface marquée associée à $X(H)$ telle que $\theta_1([X(H),\f_{1/2}])=1/2$.

\subsection{Générateurs du groupe modulaire}\label{generateurs}
Jetons un \oe il à quelques éléments du groupe des homéotopies des surfaces$_{-1}$. Soient

\begin{itemize}
\item $n$ l'homéotopie associée à la reflexion selon $\g_1$ dans $[X(P),id]$, $n$ est un élément d'ordre 2
qui fixe $\g_1$ et $\g_2$, et échange $\g_3$ et $\g_4$. Vu dans notre système de coordonnées,
$n$ correspond à la symétrie par rapport au plan $\lbrace \theta_1=0 \rbrace$.
\item $u$ la reflexion par rapport à ${\g'}_3$ dans $[X(P),id]$, $u$ échange $\g_1$ et $\g_2$. L'ensemble des points
fixes de $u$ est $\lbrace l_1=l_2 \rbrace$.
\item $t=u \cdot n$ l'homéotopie associée à la rotation d'angle $\pi/2$ en $w_1$ dans $[X(P),id]$
(voir figure~\ref{pentagone})~; $t$ échange $\g_1$ et $\g_2$, et $\g_3$ et $\g_4$.
\item $v$ l'homéotopie associée à la reflexion selon ${\g'}_1$ dans $[X(H),\f_{1/2}]$, $v$ échange $\g_2$
et $\g_3$ et son lieu des points fixes est $\lbrace l_2=l_3 \rbrace$.
\item $s=v \cdot u$ l'homéotopie qui vue dans $[X(H),\f_{1/2}]$ est la translation d'odre 6 le
long du bord de $\T_{X(H)}$, $s$ agit transitivement sur $\g_1$, $\g_2$ et $\g_3$, 
plus précisément $s(\g_1)=\g_3$ et $s(\g_3)=\g_2$.
\end{itemize}
\medskip

Nous allons montré que les éléments $n$, $s$ et $t$ engendrent le groupe modulaire $Mod_3^-$, et que ce dernier est isomorphe à $\rm{PGL(2,\Z)}$






\subsection{Domaine fondamental}

 Soit $D$ l'ensemble des points $M=[X,\f]\in\t_3^-$ tels que~:
$$0 \leq \theta_1(M) \leq 1/2\quad {\rm{et}} \quad l_1(M)\leq l_2(M).$$
L'objet de ce paragraphe est de montrer que $D$ est un domaine fondamental de $\t_3^-$ pour l'action du
groupe modulaire.

\begin{pro}\label{systole interne}
Soit $M=[X,\f]$ un point de $D$, alors $\g_1$ est une systole interne de $X$.
\end{pro}

\begin{figure}[h]
\centering
\psfrag{g}{$\gamma$}
\psfrag{g1}{$\gamma_1$}
\psfrag{d1}{${\gamma'}_1$}
\psfrag{d}{$\gamma'$}
\psfrag{w3}{$w_3$}
\psfrag{w1}{$w_1$}
\psfrag{w2}{$w_2$}

\includegraphics[width=6cm,keepaspectratio=true]{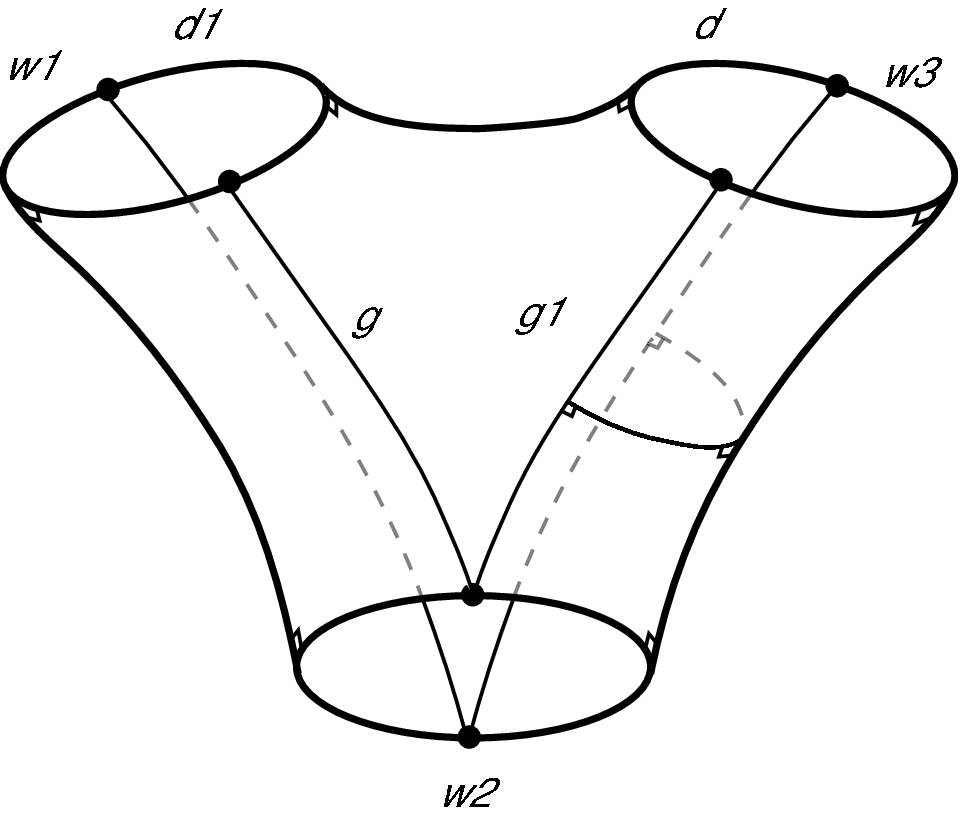}
\caption{}
\label{perpendiculaire}

\end{figure}

\begin{dem}
Soit $M=[X,\f]$ un point du plan $\lbrace \theta_1=0 \rbrace$. La surface $X$ admet un groupe d'automorphisme
$G\simeq \Z_2\times\Z_2$
avec $(r,s^+,s^-)_G=(1,1,1)$, et les systoles éventuelles de $\T_X$ sont les deux géodésiques de points fixes,
c'est-à-dire $\g_1$ et $\g_2$. Il en découle que $\g_1$ est une systole interne de tout élément $M\in D$
avec $\theta_1(M)=0.$\\
\indent Soit $M=[X,\f]$ un point de $D\setminus \lbrace \theta_1=0 \rbrace$
 tel que $X$ possède deux systoles internes~: $\g_1$ et une autre systoles
$\g$. Les géodésiques duales de $\g_1$ et $\g$ sont disjointes, et nous pouvons trouver une troisième géodésique
disjointe des deux premières~; découpons ces trois géodésiques, le pantalon obtenu a deux bords égaux. Sur ce
pantalon nous repérons sans peine $\g_1$ et $\g$ mais aussi la prependiculaire commune à $\g_1$ et $\g_X$
(voir figure~\ref{perpendiculaire}).
En particulier cette dernière n'intersecte pas $\g$.\\
\indent Maintenant découpons $\g_1$ dans $\T_X$. La géodésique $\g$ passe par $w_3$, elle traverse ensuite un
pentagone, elle ne peut en sortir par les côtés $a$ (car $\theta_1\neq 0$), $b$, $c$ ou $d$, donc elle sort
du pentagone par le point de Weierstrass (nécessairement $w_1$) du côté $e$.
En posant nos yeux sur la figure il apparaît que $\g=\g_2$.\\
\indent Comme l'ensemble $D\setminus \lbrace l_1=l_2 \rbrace$ est connexe nous concluons que $\g_1$ 
est l'unique systole interne de chacun des points de cet ensemble.
\end{dem} \hfill$\square$

\begin{cor}
Soit $M$ un point de $\t_3^-$ tel que $l_1(M) < l_{t_{\g_1}^k(\g_2)}(M)$ pour tout $k\in\Z$, alors
$\g_1$ est l'unique systole interne de $M$.
\end{cor}

\begin{dem}
A l'aide de $t_{\g_1}$ et $n$ nous revenons au cas où $M\in D$. 
\end{dem} \hfill $\square$

\begin{thm}
L'ensemble $D$ est un domaine fondamental pour l'action du groupe modulaire
$Mod_3^-$ sur l'espace de Teichmüller $\t_3^-$.
\end{thm}

\begin{dem} $D$ contient un domaine fondamental. En effet, quitte à user de $t_{\g_1}$ et $t$
nous pouvons envoyer n'importe quel point  de $\t_3^-$ sur un point $M$ vérifiant
$l_1(M) < l_{t_{\g_1}^k(\g_2)}(M)$ pour tout $k\in\Z$~; ensuite les actions de $t_{\g_1}$
et $n$ ramènent $M$ dans $D$.\\
\indent Soit $g$ un élément du groupe modulaire tel que $\ring{D}\cap (g\cdot \ring{D})\neq\emptyset$,
\begin{itemize}
\item ou bien $g\cdot \partial D =\partial D$, dans ce cas nécessairement $g(\mathcal{D}_i)=\mathcal{D}_i$
($i=1,2$) et $g=id$,
\item ou bien il existe $M\in\partial D$ tel que $g(M)\in \ring{D}$, mais ceci est impossible car les points de $D$
fixés par un élément de $Mod_3^-$ sont exactement les points frontière.
\end{itemize}
Donc $g=id$, et $D$ est bien un domaine fondamental pour l'action de groupe modulaire.
\end{dem}  \hfill $\square$ \medskip

\begin{figure}[h]
\centering
\psfrag{l}{$l_1$}
\psfrag{t}{$\theta_1$}
\psfrag{0}{$0$}
\psfrag{0,5}{$0,5$}
\psfrag{a}{$\lbrace l_1=l_3 \rbrace$}
\psfrag{b}{$\lbrace l_3=l_4 \rbrace$}
\psfrag{c}{$\lbrace l_1=l_2 \rbrace$}
\psfrag{d}{$\lbrace l_2=l_3 \rbrace$}
\psfrag{e}{$B$}
\psfrag{f}{$A$}
\psfrag{D}{$D$}

\includegraphics[width=8cm,keepaspectratio=true]{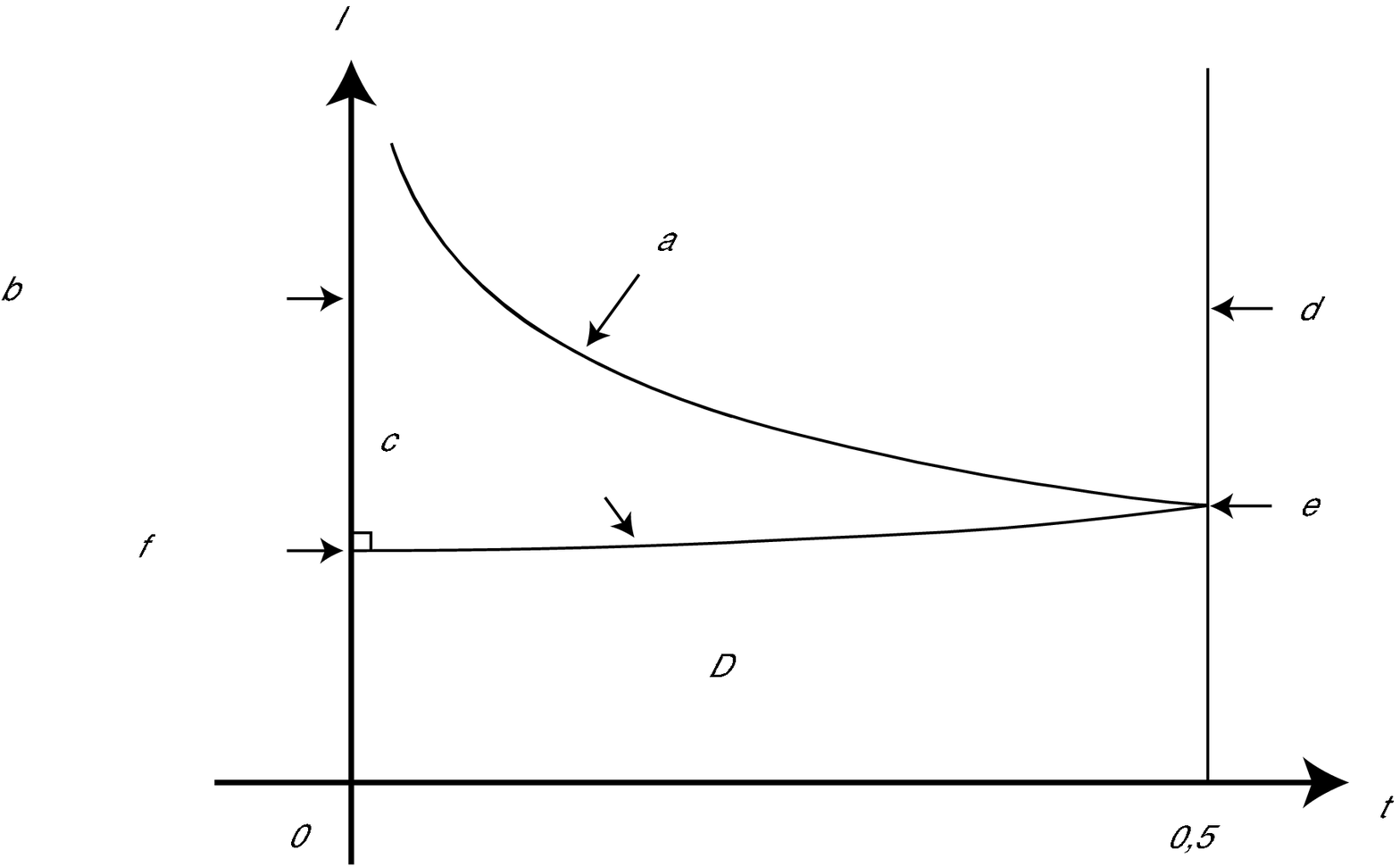}
\caption{}
\label{domaine}

\end{figure}

 Nous avons représenté en figure~\ref{domaine} le plan horizontal $\lbrace l(\g_X)=2 \rbrace$
(les échelles en abcisse et ordonné sont différentes).

Comme pour $l_1$ fixé, $l_2$ croit avec le twist $\theta_1$ (sur $\lbrace \theta_1\geq 0 \rbrace$), on observe
que les points de $\lbrace l_1=l_2,\ \theta_1\geq 0 \rbrace$ ont une ordonné qui est fonction croissante
de leur abcisse. En fait, la longueur $l_2$ s'exprime en fonction de $\theta_1$ et $l_1$ de la façon suivante~:
$$\cosh^2(l_2/2)=\cosh^2(\theta_1 l_1 /2) \frac{\cosh^2(l(\g_X)/2)+\cosh^2(l_1/2)-1}{\cosh^2(l_1/2)-1}.$$

Les points $A$ et $B$ sont respectivement les points fixes de $t$ et $s$, les surfaces associées 
admettent un groupe d'isométrie isomorphe à $\D_4$ et $\D_6$.

\begin{figure}[h]
\centering
\psfrag{l}{$l_1$}
\psfrag{t1}{$\theta_1$}
\psfrag{0}{$0$}
\psfrag{0,5}{$0,5$}
\psfrag{1}{$1$}
\psfrag{-0,5}{$-0,5$}
\psfrag{n}{$n$}
\psfrag{v}{$v$}
\psfrag{u}{$u$}
\psfrag{s}{$s$}
\psfrag{t}{$t$}
\psfrag{D}{$D$}
\psfrag{sD}{$sD$}
\psfrag{nD}{$nD$}
\psfrag{tD}{$tD$}
\psfrag{tnD}{$tnD$}
\psfrag{stnD}{$stnD$}
\psfrag{s2D}{$s^2D$}
\psfrag{ns2D}{$ns^2D$}
\psfrag{s2tnD}{$s^2tnD$}

\includegraphics[width=10cm,keepaspectratio=true]{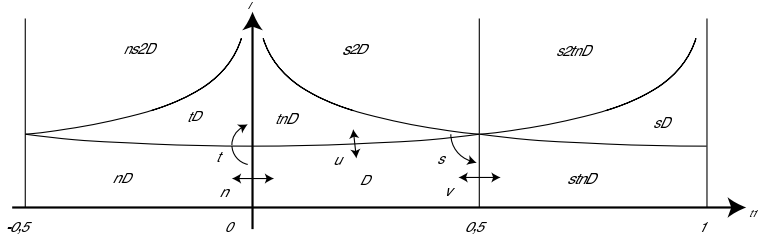}
\caption{Action des générateurs sur $D$}
\label{domaine3}

\end{figure}

Avec $n$, $t$ et $s$ nous envoyons $D$ sur n'importe lequel de ses voisins (voir figure~\ref{domaine3}),
nous en déduisons que ce sont bien des générateurs de $Mod_3^-$. 

 Nous pouvons même en dire plus. Il est clair que d'un point de vue combinatoire, le pavage obtenu par l'action de $Mod_3^-$ sur $\t_3^-$ est le même que celui obtenu par l'action classique de $\rm{PGL(2,\Z)}$ sur le demi-plan de Poincaré $\H$ ($\rm{PGL(2,\Z)}$ est le groupe d'un triangle d'angle 0, $\pi/2$, $\pi/3$). Nous en déduisons que $Mod_3^-\simeq \rm{PGL(2,\Z)}$.\par
En fait, ce résultat était déjà connu puisque
dans \cite{birman} (et \cite{birman2}) Birman et Chillingworth fournirent la présentation suivante
du groupe des homéotopies d'une surface$_{-1}$~:
$$\langle a,b,y\ |\ aba=bab,\ yay^{-1}=a^{-1},\ yby^{-1}=b^{-1},\ y^2=1,\ (aba)^4=1 \rangle.$$
Le groupe $\rm{GL(2,\Z)}$ admet pour sa part la présentation
$$\langle s,t,n\ |\ s^3t^2=1,\ ntn^{-1}t=1,\ nstn^{-1}st=1,\ n^2=1,\ t^4=1,\ s^6=1 \rangle.$$
Il y a équivalence entre ces deux présentations, il suffit de poser $n=y$, $t=aba$ et $s=ab.$

\begin{pro}[Birman, Chillingworth]
Soit $X$ une surface compacte non orientable sans bord de caractéristique -1, alors
$\rm{Homeot(X)}\simeq \rm{GL(2,\Z)}.$ \end{pro}

\begin{cor}
Le groupe modulaire $Mod_3^-$ est isomorphe à $\rm{PGL(2,\Z)}.$
\end{cor}

\begin{dem}
Hormis l'identité, seule l'involution hyperelliptique agit trivialement sur le Teichmüller, comme elle appartient
au centre du groupe modulaire $Mod_3^-$ elle s'identifie nécessairement à $-Id_2$, d'où le résultat.
\end{dem} \hfill $\square$\medskip

\subsection{Décomposition cellulaire}
\begin{figure}[Ph]
\centering
\psfrag{x}{$x$}
\psfrag{z}{$z$}
\psfrag{d4}{$\lbrace l_1=l_2, \theta_1=0  \rbrace\ $}
\psfrag{d6}{$\lbrace l_1=l_2, \theta_1=1/2\rbrace\ $}
\psfrag{lxl1}{$\lbrace l_1=l_{\g_X} \rbrace$}
\psfrag{lxd1}{$\lbrace l_{\g_X}=l_{\g'_1} \rbrace$}
\psfrag{l1d1}{$\lbrace l_1=l_{\g'_1} \rbrace$}

\includegraphics[width=6cm,keepaspectratio=true]{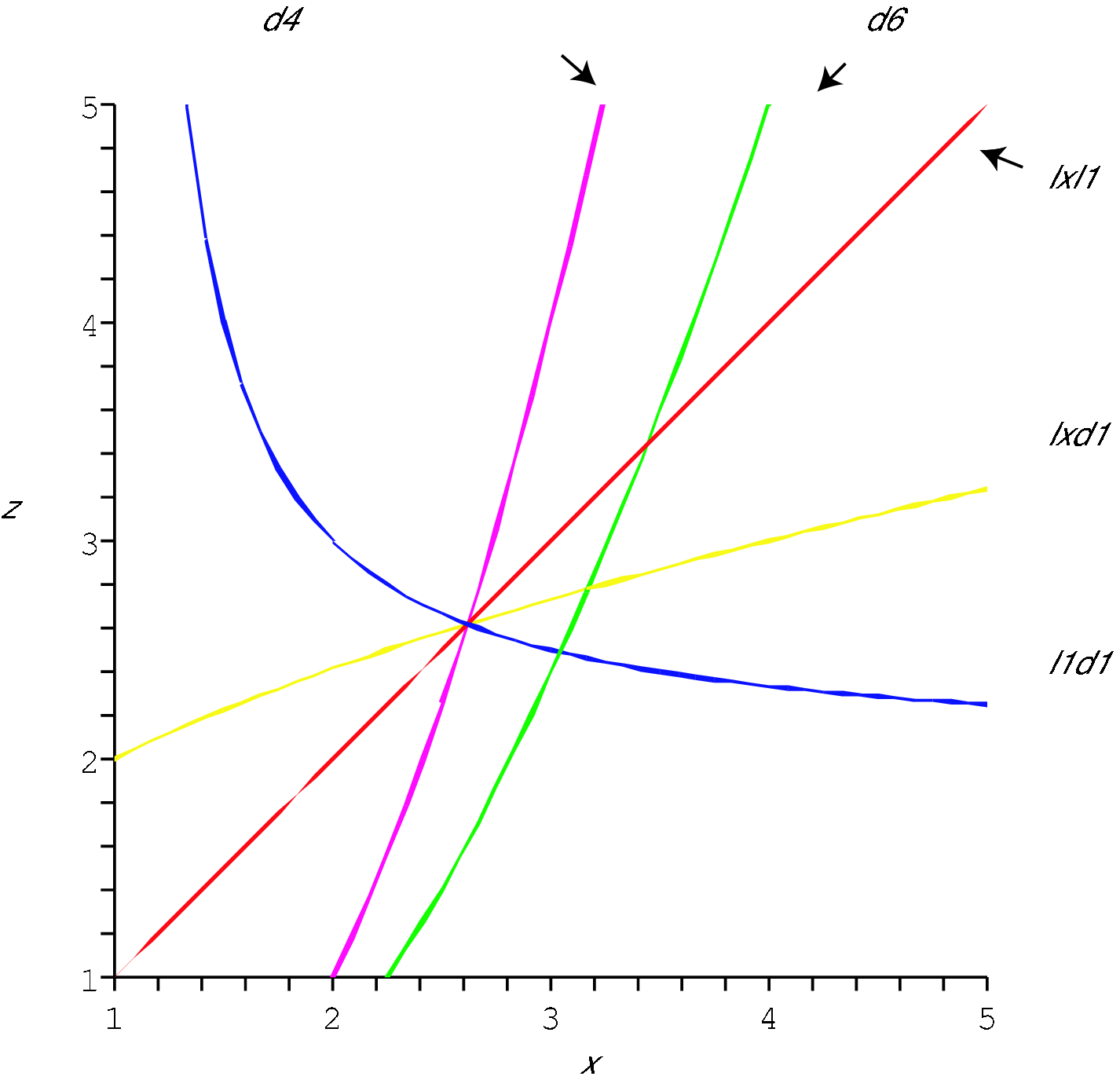}
\caption{Lieu d'egalité de longueurs}
\label{egalite}
\vspace{0.5cm}

\psfrag{t1}{$\theta_1$}
\psfrag{z}{$z$}
\psfrag{x}{$x$}
\psfrag{C1}{$C1$}
\psfrag{C2}{$C2$}
\psfrag{C3}{$C3$}
\psfrag{C4}{$C4$}
\psfrag{C5}{$C5$}
\psfrag{C6}{$C6$}

\psfrag{lxl1}{$F_{13}\subset\lbrace  l_1=l_{\g_X} \rbrace$}
\psfrag{lxd1}{$F_{23}\subset\lbrace l(\g_X)=l(\g'_1) \rbrace$}
\psfrag{l1d1}{$F_{12}\subset\lbrace l_1=l(\g'_1) \rbrace$}
\psfrag{l1l2}{$\lbrace l_1=l_2 \rbrace$}

\includegraphics[width=12cm,keepaspectratio=true]{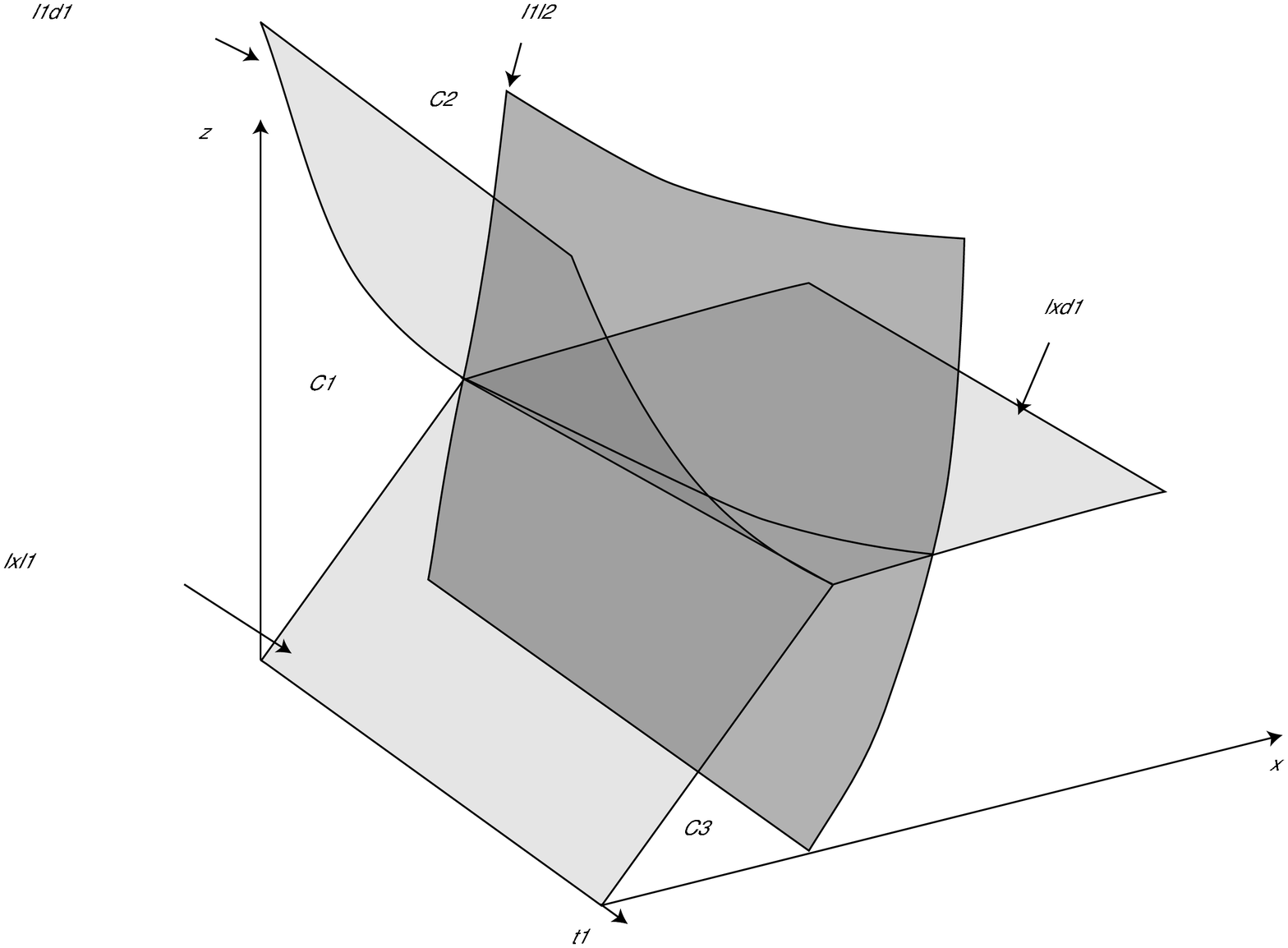}
\caption{Décomposition cellulaire}
\label{decomposition}
\end{figure}

Soit $[X,\f]$ une surface marquée de $D$, par la Proposition~\ref{systole interne} l'une des trois géodésiques
$\g_1$, ${\g'}_1$, $\g_X$ est une systole. Commençons par regarder les lieus d'égalité de longueurs de ces géodésiques.
Celles-ci ne dépendent pas du twist $\theta_1$, posons $x=\cosh^2(\g_1/2)$ et $z=\cosh^2(\g_X/2)$, en usant de l'égalité 
$\sinh({\g'}_1/2)\sinh(\g_X/2)=\cosh(\g_1/2),$
il vient~:
$$\begin{array}{clrcl}
l_{\g_1}=l_{\g_X}  & \Leftrightarrow & z       & = & x  \\
l_{\g'_1}=l_{\g_X}  & \Leftrightarrow & (z-1)^2 & = & x  \\
l_{\g_1}=l_{\g'_1} & \Leftrightarrow & z       & = & 1+\frac{x}{x-1} .\\
\end{array}$$
Par ailleurs, pour nous avons
$$\begin{array}{lclrcl}
\rm{pour}\ \theta_1=0,  & l_{\g_1}=l_{\g_2}  & \Leftrightarrow & z & = & (x-1)^2 ,  \\
\rm{pour}\ \theta_1=1/2 & l_{\g_1}=l_{\g_2}  & \Leftrightarrow & z & = & 2x^{3/2}-3x+1.  \\
\end{array}$$
Nous avons représenté en figure~\ref{egalite} les courbes associées à ces égalités, elles donnent une idée de
la configuration des lieus d'égalités de longueurs dans une coupe à twist constant de $D$ (plan vertical).

 Nous en tirons une décomposition cellulaire de $D$, une esquisse relativement fidèle
 en est donnée en figure~\ref{decomposition}.

Il y a 3 cellules de dimension 3~:
\begin{itemize}
\item $C_1=\lbrace [X,\f]\in D\ |\ l_{\g_X},l_{\g'_1}>l_{\g_1} \rbrace$, 
\item $C_2=\lbrace [X,\f]\in D\ |\ l_{\g_X},l_{\g_1}>l_{\g'_1} \rbrace$,
\item $C_3=\lbrace [X,\f]\in D\ |\ l_{\g'_1},l_{\g_1}>l_{\g_X} \rbrace$.
\end{itemize}
En excluant les faces et arêtes contenues dans les plans $\lbrace \theta_1=0 \rbrace$ et $\lbrace \theta_1=1/2 \rbrace$,
on compte 5 faces et 2 arêtes. Nous désignerons une face par un $F$ en mettant en indice le numéro de toutes
les cellules de dimension 3 adjacentes, nous procéderons de façon similaire avec les arêtes mais avec un $A$.
Pour les faces et les arêtes contenues dans les plans $\lbrace \theta_1=0 \rbrace$ et $\lbrace \theta_1=1/2 \rbrace$,
nous adopterons le même système mais avec des $F^{0}$, des $A^{0}$, des $F^{1/2}$ et des $A^{1/2}$.
Tous les points d'une cellule ont même(s) systole(s), les tableaux ci-dessous donne la (ou les) systole(s)
associée(s) à une cellule.
\begin{center}
\begin{tabular}{|l|c|l|c|l|c|}\hline
cellule    & systole     & face      & systoles      & arête        & systoles           \\ \hline
$C_1$      & $\g_1$      & $F_{12}$  & $\g_1,\g'_1$  & $A_{123}$    & $\g_1,\g'_1,\g_X$  \\ \hline
$C_2$      & $\g'_1$     & $F_{2} $  & $\g'_1,\g'_2$ & $A_{23}$     & $\g'_1,\g'_2,\g_X$ \\ \hline
           &             & $F_{23}$  & $\g'_1,\g_X$  &              &                    \\ \hline
$C_3$      & $\g_X$      & $F_{3} $  & $\g_X$        &              &                    \\ \hline
           &             & $F_{13}$  & $\g_1,\g_X$   &              &                    \\ \hline
\end{tabular}\end{center}
\medskip

\begin{center}
\begin{tabular}{|l|c|l|c|l|c|}\hline
 face      &  systole      & arête         & systoles      & sommet      & systoles                      \\ \hline
$F_1^{0}$  & $\g_1$        & $A^{0}_{12}$  & $\g_1,\g'_1$  & $[X(P),id]$ & $\g_1,\g'_1,\g_2,\g'_2,\g_X$  \\ \hline
$F_2^{0}$  & $\g'_1$       & $A^{0}_{2}$   & $\g'_1,\g'_2$ &             &                               \\ \hline
$F_3^{0}$  & $\g_X$        & $A^{0}_{3}$   & $\g_X$        &             &                               \\ \hline
           &               & $A^{0}_{13}$  & $\g_1,\g_X$   &             &                               \\ \hline
\end{tabular}\end{center}
\medskip

\begin{center}
\begin{tabular}{|l|c|l|c|l|c|}\hline
 face        &  systole & arête         & systoles      & sommet      & systoles                      \\ \hline
$F_1^{1/2}$  & $\g_1$   & $A^{1/2}_{12}$& $\g_1,\g'_1$  & $S_{123}$     & $\g_1,\g'_1,\g_X$             \\ \hline
$F_2^{1/2}$  & $\g'_1$  & $A^{1/2}_{2}$ & $\g'_1,\g'_2,\g'_3$ &       &            \\ \hline
             &          & $A^{1/2}_{23}$& $\g'_1,\g_X$  & $[X(H),\f_{1/2}]$ & $\g'_1,\g'_2,\g'_3,\g_X$ \\ \hline
$F_3^{1/2}$  & $\g_X$   & $A^{1/2}_{3}$ & $\g_X$        &             &             \\ \hline
             &          & $A^{1/2}_{13}$& $\g'_1,\g_X$  &             &             \\ \hline
\end{tabular}
\end{center}




\subsection{Points eutactiques, parfaits et extrêmes}
 Nous laissons au lecteur le soin de lire les articles \cite{bavard1} et \cite{bavard3} de C. Bavard afin de prendre connaissance des différentes définitions et résultats classiques de cette théorie. Donnons cependant un analogue du théorème de Voronoï, il est énoncé pour le cas orientable dans les articles cités mais reste valable dans le cadre des surfaces hyperboliques compactes de caractéristiques -1~:

\begin{thm}[Bavard] Une surface hyperbolique fermée de caractéristique -1 est extrême si est seulement si elle est parfaite et eutactique.
\end{thm}
  
 D'après un résultat de H. Akrout (\cite{akrout}), dans le cas orientable la systole est une fonction de Morse topologique sur l'espace de Teichmüller~; ses points critiques sont exactement les points eutactiques et par conséquent, sont intéressants à classer en eux-mêmes. Ceci reste vrai dans la cas non orientable, aussi nous recherchons tous les points eutactiques.

\begin{pro}
Il y a exactement deux surfaces$_{-1}$ eutactiques~: $X(P)$ et $X(H)$.
\end{pro}

\begin{dem} Il suffit de déterminer les points eutactiques dans le domaine $D$.
Regardons les gradients des fonctions longueur des systoles éventuelles de
surface$_{-1}$, nous voyons directement
$$\nabla l_1 =  \left( \begin{array}{c} 0 \\ 1 \\ 0  \end{array} \right),
\nabla l_{\g_X} =\left( \begin{array}{c} 0 \\ 0 \\ 1  \end{array} \right),
\nabla l_{\g'_1}=\left( \begin{array}{c} 0 \\ >0 \\ <0  \end{array} \right).$$
 Ensuite en dérivant par rapport à $l_1$ l'égalité $$\cosh(l_2/2)=\cosh(\theta_1 l_1/2)\frac{\cosh(l_X/2)}{\sinh(l_1/2)},$$
nous trouvons
\begin{eqnarray*}
\frac{\partial l_2}{2\partial l_1}\sinh(l_2/2) & = & \frac{\theta_1}{2} \sinh(\theta_1 l_1/2)\frac{\cosh(l_{\g_X}/2)}{\sinh(l_1/2)}
-\cosh(\theta_1l_1/2)\cosh(l{\g_X}/2)\frac{\cosh(l_1/2)}{2\sinh^2(l_1/2)} \\
              & \leq & \frac{\cosh(l_{\g_X}/2)\cosh(\theta_1l_1/2)}{2} \left[ \frac{1}{2\sinh(l_1/2)}
-\frac{\cosh(l_1/2)}{\sinh^2(l_1/2)} \right] \\
 & < & 0.
\end{eqnarray*}
En dérivant par rapport aux autres variables il vient 
$$\left\lbrace
\begin{array}{lc}
\nabla l_{\g_2}([X,\f])=\left( \begin{array}{c} 0 \\ <0 \\ >0  \end{array} \right)\ \rm{si} & \theta_1([X,\f])=0 \\
\nabla l_{\g_2}([X,\f])=\left( \begin{array}{c} >0 \\ <0 \\ >0  \end{array} \right)\ \rm{si} & 1/2\geq\theta_1([X,\f])>0
\end{array}.
\right.$$
Enfin via l'égalité de la remarque~\ref{rem dualite} nous trouvons
$$\left\lbrace
\begin{array}{lc}
\nabla l_{\g'_2}([X,\f])=\left( \begin{array}{c} 0 \\ <0 \\ \ast  \end{array} \right)\ \rm{si} & \theta_1([X,\f])=0 \\
\nabla l_{\g'_2}([X,\f])=\left( \begin{array}{c} >0 \\ <0 \\ \ast  \end{array} \right)\ \rm{si} & 1/2\geq\theta_1([X,\f])>0
\end{array}.
\right.$$

Nous allons étudier les différentes configuration de systoles pour un point $[X,\f]$ de $D$.
\begin{itemize}
\item Si $X$ possède exactement une systole, alors le gradient de la fonction longueur associée est non nul
et $X$ n'est pas une surface eutactique.
\item Si $X$ possède exactement deux systoles, et si celles-ci sont $\g_1$ et $\g'_1$, ou $\g_1$ et $\g_X$, ou
$\g'_1$ et $\g_X$, alors les vecteurs gradient des fonctions longueur des systoles ne sont visiblement
pas colinéaires et la surface n'est pas eutactique. Si $X$ contient exactement deux systoles qui sont
$\g'_1$ et $\g'_2$, alors en regardant les gradients des fonctions longueur dans le système de coordonnées
$(l_{\g'_1},l_{\g'_2},l_{\g'_3})$ on s'aperçoit qu'ils ne sont pas colinéaires, et par suite la surface pas eutactique.
\item Si $X$ possède exactement trois systoles, alors ou bien les systoles de $X$ sont $\g'_1$, $\g'_2$ et $\g'_3$, et
en se plaçant dans le système de coordonnées $(l_{\g'_1},l_{\g'_2},l_{\g'_3})$ il est clair que $X$ n'est pas
eutactique~; ou bien les systoles de $X$ sont $\g_1$, $\g'_1$ et $\g_X$, et avec les coordonnées usuelles on voit
que $X$ n'est pas eutactique~; ou bien les systoles de $X$ sont $\g'_1$, $\g'_2$ et $\g_X$, dans ce cas
$\theta_1([X,\f])>0$ et en regardant les coordonnées des vecteurs gradient des fonctions longueurs on voit encore une fois
que $X$ n'est pas eutactique.
\item Si $X$ possède exactement quatre systoles, alors $X=X(H)$. Comme $X(H)$ est un point fixe isolé d'un sous-groupe de $Mod_3^-$ c'est un point eutactique (voir \cite{bavard3}).
\item Si $X$ possède exactement cinq systoles, alors $X=X(P)$. Comme $X(P)$ est elle aussi un point fixe isolé d'un sous-groupe du groupe modulaire, elle est eutactique.


\end{itemize}

\end{dem}
 
Avec des calculs triviaux nous déterminons aussi les surfaces parfaites et ainsi les surfaces extrêmes de $\t_3^-$.

\begin{pro}
Il y a une et une seule surfaces$_{-1}$ parfaite~: $X(H)$.
\end{pro}

\begin{cor}
Il y a une et une seule surfaces$_{-1}$ extrême~: $X(H)$.\par
 Le maximum global de la fonction systole sur $\t_3^-$ (modulo l'action du groupe modulaire)
est atteint en $[X(H),\f_{1/2}]$ et $$\cosh(sys(X(H)))=\frac{5+\sqrt{17}}{2}.$$
\end{cor}

\begin{rem}
Une fois connus tous les points eutactiques de $\t_3^-$ via la formule donnée dans l'article \cite{akrout}, nous pouvons calculer la caractéristique d'Euler de $Mod_3^-$.
\end{rem}

\section{Cas des surfaces hyperboliques de caractéristique~-1 à bords}

De notre étude de $\t_3^-$ nous allons déduire une description des espaces
de Teichmüller du tore à un bord, du plan projectif à deux bords et de la bouteille de Klein à un bord
(les longueurs des bords seront supposées fixées). \medskip

 Avant toute chose, pour toutes les surfaces compactes possédant $k$ bords, nous supposerons déjà fixées
une même dénomination $(b_i)_{i=1\ldots k}$ des bords, ainsi que leurs longueurs $l(b_i)_{i=1\ldots k}$.
Sauf mention du contraire, un homéomorphisme entre deux surfaces à bords préservera les dénominations des bords.

\subsection{Espaces de Teichmüller des surfaces à bords}
 Nous procédons comme nous l'avions fait dans le cas des surfaces sans bord.
Nous fixons une surface compacte $R$ orientable ou non, de genre $g$ avec $k$ bords
(les longueurs des bords, nous venons de le dire, sont fixées).
Une surface marquée est alors un couple $(S,\f)$, où $S$ est une surface et $\f:R\mapsto S$ un homéomorphisme,
préservant l'orientation si $R$ est orientable.
Deux surfaces marquées $(S,\f)$ et $(S',\f')$ seront équivalentes lorsque
l'homéomorphisme $\f'\circ\f^{-1}$ sera isotope à une isométrie.
On définit l'espace $\t_{g,k}$ (resp. $\t_{g,k}^-$) des surfaces compactes orientables (resp. non orientables)
de genre $g$ à $k$ bords de longueurs fixées, comme l'ensemble des classes
d'équivalence de surfaces marquées. Cet espace possède une structure de variété $C^\infty$
déterminée par les fonctions longueurs.\medskip

 Nous allons montrer que par auto-recollement des bords nous pouvons plonger ces espaces de Teichmüller
de surfaces à bords ($\t_{1,1}$, $\t_{1,2}^-$ et $\t_{2,1}^-$ par exemple) dans les espaces de Teichmüller
de surfaces sans bords (ici $\t_3^-$).\medskip

 Soit $S$ une surface compacte avec $k$ bords, notons $\bar{S}$ la surface obtenue par auto-recollement
des bords. Nous avons le lemme suivant~: 

\begin{lem}\label{lem_1}
Soit $\f$ un homéomorphisme de $S$,
$\f$ induit un homéomorphisme $\bar{\f}$ de $\bar{S}$ unique à isotopie près.
De plus, si $\psi$ est un homéomorphisme de $S$ isotope à $\f$, alors $\bar{\psi}$ et $\bar{\f}$ sont
isotopes.
\end{lem}

 Afin de donner une esquisse de démonstration de ce lemme nous aurons besoin d'un résultat bien connu.

\begin{pro}\label{isotopie}
Un homéomorphisme du cercle est isotope à une isométrie.
\end{pro}


\noindent\textbf{Esquisse de démonstration du lemme~\ref{lem_1}}
 La restriction de $\f$ à chacun des bords est isotope à un homéomorphisme commutant avec
l'auto-recollement de ce bord (par la proposition précédente), ainsi $\f$ induit un homéomorphisme
de $S$ commutant avec l'auto-recollement des $b_i$ ($i=1\ldots k$),
donc induit un homéomorphisme $\bar{\f}$ de $\bar{S}$.
 Avec un peu de topologie, il apparaît que si $\f$ est isotope à $\psi$, un autre homéomorphisme de $S$,
alors $\bar{\f}$ et $\bar{\psi}$ sont isotopes. Nous concluons que l'homéomorphisme induit est bien unique
à isotopie près. 
\hfill$\square$
\medskip

 Ainsi, à une surface marquée $[S,\f]$ correspond une unique surface marquée $[\bar{S},\bar{\f}]$.
Autrement dit, nous définissons (modulo l'action des groupes modulaires) une unique application
$$\begin{array}{cccc}
\mathcal{R}: & \t_{g,k}\ \rm{(resp. \t_{g,k}^-)}  & \longrightarrow &  \t_{2g+k}^-\ \rm{(resp. \t_{g+k}^-)}\\
             &  [S,\f]    & \longmapsto & [\bar{S},\bar{\f}] \\ 
\end{array}.$$
Notons $f_1$,...,$f_k$ les géodésiques issues de l'auto-recollement des bords, 
$\mathcal{R}$ est à valeurs dans l'ensemble $\mathcal{F}=\lbrace l(f_i)=l(b_i)/2,\ i=1\ldots k \rbrace$.\medskip

 Inversement, soit $S$ une surface compacte non orientable sans bord, et soient $f_1,\ldots,$ $f_k$ des
géodésiques de $S$ non orientables et disjointes.
Si $\hat{S}$ désigne la surface obtenue par découpe de ces $k$ géodésiques, alors le lemme suivant vient directement.

\begin{lem}\label{lem_2}
Soit $\f$ un homéomorphisme de $S$ fixant (à homotopie près) les courbes $(f_i)_{i=1\ldots k}$,
$\f$ induit un homéomorphisme $\hat{\f}$ de $\hat{S}$ unique à isotopie près.
De plus, si $\psi$ est un homéomorphisme de $S$ isotope à $\f$, alors $\hat{\psi}$ et $\hat{\f}$ sont
isotopes.
\end{lem}

 Nous pouvons alors définir une application
$$\begin{array}{cccc}
\mathcal{D}: & \mathcal{F}=\lbrace l(f_i)= l(b_i)/2,\ (i=1\ldots k \rbrace \subset \t_g^- & \longrightarrow &  \t_{(g-k)/2,k}\ \rm{ou}\ \t_{g-k,k}^- \\
             &  [S,\f]    & \longmapsto & [\hat{S},\hat{\f}] \\ 
\end{array},$$
l'espace d'arrivée dépendant de la nature topologique des géodésiques.\medskip

 Soient $\t1$ un espace de Teichmüller de surfaces à bords, et $\t2$ l'espace de Teichmüller
des surfaces sans bords de même caractéristique. Les applications
$\mathcal{D}$ et $\mathcal{R}$ existantes entre ces deux espaces vérifient $\mathcal{D}\circ\mathcal{R}=id_{\t_1}$ et
$\mathcal{R}\circ\mathcal{D}=id_{\mathcal{F}}$. En particulier l'image de $\mathcal{R}$ est exactement l'ensemble 
$\mathcal{F}$ des points de $\t2$ vérifiant les conditions de longueur. D'une manière générale les lieus d'égalité
de longueurs ne constituent pas des sous-variété du Teichmüller, néanmoins dans le cas présent les géodésiques
$(f_i)_{i=1\ldots k}$ sont disjointes, et par un petit lemme nous concluons que $\mathcal{F}$ est une sous-variété
de $\t_g^-$.

\begin{lem} 
Soit $R$ une surface hyperbolique compacte, soit $(f_i)_{i=1\ldots k}$ une famille
de géodésiques disjointes de $R$. L'ensemble
$\mathcal{F}=\lbrace l(f_i)=l(b_i)/2,\ i=1\ldots k \rbrace$ est une sous-variété lisse de codimension $k$
de l'espace de Teichmüller de $R$. 
\end{lem}

\begin{dem}
Nous pouvons compléter cette famille en une famille de géodésiques
disjointes, telle que les coordonnées de longueur et éventuellement de twist associées établissent
un difféomorphisme entre $\t(R)$ et un produit $(\R_+^\ast)^{\alpha}\times \R^\beta$
avec $\alpha+\beta=dim(\t(R))$. L'assertion s'en déduit immédiatement.\nolinebreak$\square$
\end{dem}\medskip

 En choisissant de bons systèmes de coordonnées de $\t1$ et $\t2$, il est désormais clair que les applications
de découpe et de recollement sont de classe $C^\infty$. Par conséquent,

\begin{pro}
Les applications $\mathcal{R}$ et $\mathcal{D}$ établissent des difféomorphismes entre $\t1$ et
$\mathcal{F}\subset\t2$. En particulier, $R$ est un plongement de $\t1$ dans $\t2$.
\end{pro}

 En ce qui concerne le groupe des homéotopies d'une surface à bord $S$,
par les lemmes \ref{lem_1} et \ref{lem_2} il est isomorphe au groupe $Fix((f_i)_{i=1\ldots k})$ égal à
$$\lbrace h\in{\rm{Homeot(\bar{S})}}~;\ h\ fixe\ (\grave{a} \ homotopie\ pr\grave{e}s)\ les\ g\acute{e}od\acute{e}siques\
f_i\ i=1\ldots k\rbrace.$$
L'application $\mathcal{R}$ est
équivariante pour l'action de $\rm{Homeot(S)}$. \medskip

 Enfin, nous aurions pu définir une application d'auto-recollement de $k_1$ bords (avec $k_1<k$) d'une surface à $k$ bords,
nous retrouverions avec les mêmes preuves tous les résultats obtenus.

\subsection{Le tore à un bord}
 L'espace de Teichmüller $\t_{1,1}$, des tores à un bord de longueur
fixée $l(b_1)>0$, s'identifie au plan $\lbrace l_{\g_X}=l(b_1)/2 \rbrace \subset \t_3^-$ via le plongement $\mathcal{R}$.
Le groupe des homéotopies est isomorphe à celui des surfaces$_{-1}$ puisque tout homéomorphisme
de surface$_{-1}$ fixe l'ovale, ainsi il est isomorphe à 
$\rm{GL(2,\Z)}$, et de même $Mod_{1,1}\simeq\rm{PGL(2,\Z)}$. \medskip

 Un domaine fondamental $D_{l(b_1)}$ de l'action de $Mod_{1,1}$ sur $\t_{1,1}$ s'obtient en intersectant $D$ avec
le plan $\lbrace l_{\g_X}=l(b_1)/2 \rbrace$ (voir figure~\ref{domaine}), et l'action des générateurs se lit
directement sur la figure~\ref{domaine3}.\medskip

 Considérons la décomposition cellulaire triviale, c'est-à-dire constituée d'une face $F$, trois arêtes  
$A_1=\lbrace l_3=l_4 \rbrace$, $A_2=\lbrace l_1=l_2 \rbrace$, et $A_3=\lbrace l_2=l_3 \rbrace$, et
deux sommets $S_{12}$ et $S_{23}$. Les points de la face, de même que ceux des arêtes $A_1$ et $A_3$
admettent une seule systole, $\g_1$. Les points de l'arête $A_2$ possèdent deux systoles, $\g_1$ et $\g_2$~;
mais aucun d'eux n'est un point eutactiques car $\frac{\partial l_2}{\partial \theta_1}\neq 0$.
Enfin, les systoles de $S_1$ sont $\g_1$ et $\g_2$, et celles de $S_2$, $\g_1$, $\g_2$ et $\g_3$.
Ces deux surfaces, en tant que points fixes isolés de sous-groupes du groupe modulaire sont eutactiques.
Trivialement la seule surface parfaite correspond au sommet $S_{23}$. Nous retrouvons ainsi ce résultat
de P. Schmutz Schaller (\cite{schmutz})~:

\begin{thm}
L'espace $\t_{1,1}$ des tores à un bord de longueur fixée ($l(b_1)>0$) contient une seule surface extrême (modulo l'action du groupe modulaire).
Cette surface est modélisée sur l'hexagone régulier (figure~\ref{hexagone}), et la longueur $s$ de sa systole est donnée par
$$\cosh(s/2)=\cosh(l(b_1)/6)+1/2.$$
\end{thm}

\subsection{La bouteille de Klein à bord}
 Regardons maintenant l'espace $\t_{2,1}^-$ des bouteilles de Klein à bord de longueur fixée $l(b_1)>0$.
Via le plongement $\mathcal{R}$ et en supposant que le bord auto-recollé se confond avec la géodésique $\g'_1$, cet espace s'identifie à $\lbrace l_{\g'_1}=l(b_1)/2 \rbrace\subset \t_3^-$. Remarquons que cet ensemble est une surface réglée, son intersection avec un plan horizontal
est une droite parallèle à l'axe du twist.\par

 Les éléments de $Mod_3^-$ fixant $\g'_1$ n'agissent que sur la coordonnée de twist $\theta_1$ car ils fixent aussi $\g_X$. Sur la figure~\ref{domaine3} nous voyons aisément quels sont ces éléments, ce sont ceux envoyant $D$ sur un autre domaine fondamental bordé par le plan $\lbrace l_1=0 \rbrace$. Leur groupe est engendré par $n$ et $v=stn$, les reflexions selon les droites $\lbrace \theta_1=0\rbrace$ et $\lbrace \theta_1=1/2\rbrace$~; il se présente de la façon suivante~:
 $\langle n,v\ |\ n^2=v^2=1 \rangle$.
 
 Un domaine fondamental pour l'action de ce groupe sur $\lbrace l_{\g'_1}=l(b_1)/2 \rbrace$ est donc $\lbrace l_{\g'_1}=l(b_1)/2, 0\leq\theta_1\leq 1/2 \rbrace$.\par
 
 Une décomposition cellulaire intéressante de ce domaine, est celle constituée de deux faces 
\begin{itemize}
\item $F_1=\lbrace l_{\g'_2}>l_{\g_1} \rbrace$,
\item $F_2=\lbrace l_{\g'_2}<l_{\g_1} \rbrace$,
\end{itemize}
de cinq arêtes
 \begin{itemize}
\item $A_{11}=\lbrace l_{\g'_1}=l(b_1)/2,\ \theta_1=0,\ l_{\g'_2}>l_{\g_1}\rbrace, $
\item $A_{12}=\lbrace l_{\g'_1}=l(b_1)/2,\ \theta_1=1/2,\ l_{\g'_2}>l_{\g_1}\rbrace$,
\item  $A_{21}=\lbrace l_{\g'_1}=l(b_1)/2,\ \theta_1=0,\ l_{\g'_2}<l_{\g_1}\rbrace $,
\item $A_{22}=\lbrace l_{\g'_1}=l(b_1)/2,\ \theta_1=1/2,\ l_{\g'_2}<l_{\g_1}\rbrace$,
\item $A_{3}=\lbrace l_{\g'_1}=l(b_1)/2,\ l_{\g'_2}=l_{\g_1}\rbrace, $,
\end{itemize}
et de deux sommet $S_{1}$ (le sommet à twist nul) et $S_{2}$ (le sommet à twist 1/2).\par
Montrons que l'ensemble $A_{3}$ forme bien une arête. Nous calculons $l_2$ de deux manières~:
\begin{eqnarray*}
 \cosh^2(l_2/2) & = & \sinh^2(l_{\g'_2}/2)\sinh^2(l_X/2) \\
 \cosh^2(l_2/2) & = & \cosh^2(\theta_1 l_1/2)\coth^2(l_1/2)\coth^2(l_{\g'_1}/2) \\
                & = & \cosh^2(\theta_1 l_1/2) \frac{\sinh(l_{\g'_1}/2)\sinh(l_X/2)}{\sinh^2(l_1)}.
\end{eqnarray*}
En utilisant successivement la première puis la deuxième de ces égalités, il vient~:
\begin{eqnarray*}
\cosh^2(l_{\g'_2}/2) & = & \frac{\cosh^2(l_2/2)}{\sinh^2(l_X/2)} +1 \\
                    &=& \frac{\cosh^2(\theta_1 l_1/2) \sinh(l_{\g'_1}/2)}{\cosh^2(l_1/2)-1} +1 .
\end{eqnarray*}
Au final, nous trouvons que $l_{\g'_2}=l_{\g_1}$ équivaut à 
$$\cosh^2(l_{\g'_1}/2)   =  \left[ \frac{\sinh^2(l_1/2)}{\cosh(\theta_1 l_1/2)}\right]^2.$$
Considérons $\theta_1$ fixé, la fonction $x\mapsto \frac{\sinh^2(x)}{\cosh(\theta_1 x)}$ est strictement croissante sur $\R_+^\ast$ et tend vers 0 en 0, vers $+\infty$ en $+\infty$. Pour tout
$\theta_1$ il existe donc un unique $l_1$ tel que $l(\g'_2)=l(\g_1)$.\medskip
 
 Comme la coordonnée de twist est comprise entre 0 et 1/2, les seules systoles possibles sont $\g'_2$, $\g'_3$ et $\g_1$, car les autres géodésiques de la bouteille de Klein à bord se déduisent de celles-ci par le twist de Dehn $t_{\g_1}$. Soulignons aussi que si $\g'_3$ est une systole, alors $\theta_1=1/2$ (sinon $l(\g'_2)<l(\g'_3)$). Nous pouvons maintenant donner la systole de chaque composante~:
\begin{center}
\begin{tabular}{|l|c|l|c|l|c|}\hline
face       & systole     & arête     & systoles      & sommet       & systoles           \\ \hline
$F_1$      & $\g_1$      & $A_{11}$  & $\g_1$        & $S_1$        & $\g_1,\g'_2$       \\ \hline
           &             & $A_{12}$  & $\g_1$        &              &                    \\ \hline
           &             & $A_{3}$   & $\g_1,\g'_2$  & $S_3$        & $\g_1,\g'_2,\g'_3$ \\ \hline
$F_2$      & $\g'_2$     & $A_{21}$  & $\g'_2$       &              &                    \\ \hline
           &             & $A_{22}$  & $\g'_2,\g'_3$ &              &                    \\ \hline
\end{tabular}
\end{center}

\medskip

 Déterminons les points eutactiques. Les surfaces possédant une seule systole ne sont pas des points eutactiques puisque les fonctions longueurs sont strictement convexes. Les points ayant $\g'_2$ et $\g'_3$ pour systoles ne sont eux non plus pas eutactiques puisque $(l_{\g'_2},l_{\g'_3})$ est un système de coordonnées de $\lbrace l_{\g'_1}=l(b_1)/2 \rbrace$. Il ne reste plus qu'à étudier les deux sommets. Nous observons que 
$$\nabla l_1= \left( \begin{array}{c} 0 \\ 1 \end{array}\right)\ \rm{et}\ \nabla l_{\g'_2}(S_1)= \left( \begin{array}{c} 0 \\ <0 \end{array}\right),$$ donc $S_1$ est eutactique.
Penchons nous sur le cas de $S_2$, l'expression de $l_1$ est symétrique en les coordonnées $(l_{\g'_2},l_{\g'_3})$, donc le gradient $\nabla l_1(S_2)$ a même composantes en $\frac{\partial}{\partial l_{\g'_2}}$ et $\frac{\partial}{\partial l_{\g'_2}}$, il est facile de voir qu'elles sont négatives ce qui implique $S_2$ eutactique.\par
 Evidemment, $S_2$ est le seul point parfait, et aussi le seul point extrême. Pour calculer sa systole, nous découpons $\g'_2$ et $\g'_3$, nous obtenons un pantalon à deux bords égaux. Ce pantalon se divise deux hexagones droits, nous travaillons dans l'un d'eux. Notons $s$ la longueur des systoles et $z$ la longueur de la perpendiculaire commune à $\g'_2$ et $\g'_3$, nous avons (formule exprimant la distance de déplacement)
$\sinh(s/4)=\sinh(z/2)\cosh(s/2)$ , mais $$\sinh(z/2)=\cosh(l(b_1)/4)/(2\sinh(s/2)\cosh(s/2)),$$
d'où $\cosh^2(s/2)=\cosh(l(b_1)/4)/2+1$ et $\cosh(s)=\cosh(l(b_1)/4)+1$. 




\begin{pro}
Modulo l'action du groupe modulaire, l'espace $\t_{2,1}^-$ a deux points eutactiques, un point parfait et un point extrême. Soit $s$ la systole de ce point, la valeur de $s$ est donnée par
$$\cosh(s)=\cosh(l(b_1)/4)+1.$$
\end{pro}

\subsection{Le plan projectif à deux bords}

 Considérons $\t_{1,2}^-$ l'espace de Teichmüller des plans projectifs à deux bords de
longueurs $l(b_1)$ et $l(b_2)$ fixées. Par l'application d'auto-recollement $\mathcal{R}$,
cet espace se plonge dans $\t_3^-$. Munissons $\t_3^-$ du système de coordonnées
$(l_{\g'_1},l_{\g'_2},l_{\g'_3})$, et
supposons que les bords des plans projectifs sont identifiés à $\g'_1$ et $\g'_2$ après auto-recollement.
Alors, l'image de $\t_{1,2}^-$ par $\mathcal{R}$ est la droite $\lbrace l_{\g'_1}=l(b_1)/2,\ l_{\g'_2}=l(b_2)/2 \rbrace$,
en particulier $\t_{1,2}^-$ est paramétré par la longueur $l_{\g'_3}$ (c'est trivial, voir figure~\ref{plan projectif})\medskip

 Nous pouvons montrer simplement que $Mod_{1,2}^-$ s'identifie à $\langle n \rangle\simeq \Z/2\Z$. En effet, on obtient un plan projectif à deux bords en découpant une géodésique non orientable d'une bouteille de Klein à bord, donc $Mod_{1,2}^-$ peut être vu comme le sous-groupe de $Mod_{2,1}^-$ stabilisant une certaine géodésique.\par
 Ceci dit, des présentations des groupes des homéotopies des surfaces compactes non orientables furent données
par J.S. Birman et D.R.J. Chillingworth dès 1971. Plus récemment, M. Korkmaz fournit des générateurs
du groupes des homéotopies pour les surfaces non orientables pointées, et dans le cas du plan projectif moins deux
points montra que ce groupe est isomorphe à $\Z/2\Z\times \Z/2\Z$. Or, nous savons par un théorème
de D.J. Sprows, que le groupe des homéotopies d'une surface compacte privée de $k$ points est isomorphe
au groupe des homéotopies de cette surface privée de $k$ disques ouverts. Par conséquent, 

\begin{pro}
Le groupe des homéotopies du plan projectif à deux bords est isomorphe à $\Z/2\Z\times \Z/2\Z$.
\end{pro}

\begin{cor}
Le groupe $Mod_{1,2}^-$ est isomorphe à $\Z/2\Z$.
\end{cor}

 Soit $Z$ un plan projectif à deux bords, le groupe $\rm{Homeot}(Z)$ est isomorphe
au sous-groupe de $\rm{Homeot(\bar{Z})}$ formé des éléments fixant (à homotopie près) chacune 
des courbes $\g'_1$ et $\g'_2$. Nous observons que l'identité $id$,
l'involution hyperelliptique $\iota_Z$, l'homéotopie $n$ définie en \textsection~\ref{generateurs}
et le produit $\iota \cdot n$, vérifient tous cette contrainte. Par conséquent, les homéotopies de $Z$
qui s'en déduisent composent l'ensemble du groupe des homéotopies.\medskip

 Décrivons l'action de ces éléments.
L'involution hyperelliptique renverse l'orientation des géodésiques de $Z$, ainsi que celle des bords. 
L'homéotopie $n$ conserve l'orientation de $b_1$ et renverse celle
de $b_2$, elle échange les deux géodésiques fermées simples de $Z$. $n$ est ce qu'on appelle un glissement
de bord (\guillemotleft boundary slide\guillemotright\  en anglais), cela consiste à faire passer un bord d'une
surface non orientable à travers un plan projectif.
Ce concept fut introduit par Lickorish, nous renvoyons à \cite{korkmaz} pour plus de détails.
Notons qu'en général un glissement n'est pas involutif car son carré est un twist de Dehn le long d'une courbe,
cependant dans le cas présent cette courbe est un bord.

\begin{figure}[Ht]
\centering
\psfrag{o}{$\tilde{\gamma}_X$}
\psfrag{e1}{$\tilde{b}_{11}$}
\psfrag{e2}{$\tilde{b}_{21}$}
\psfrag{e12}{$\tilde{b}_{12}$}
\psfrag{e22}{$\tilde{b}_{22}$}
\psfrag{g3}{$\tilde{\g'}_3$}
\psfrag{g4}{$\tilde{\g'}_4$}
\psfrag{P}{$\mathcal{P}$}

\includegraphics[width=4cm,keepaspectratio=true]{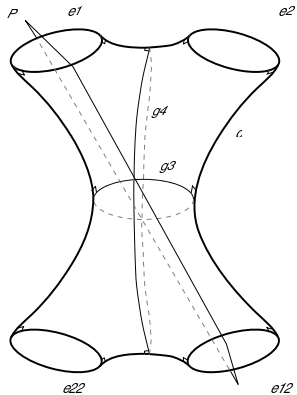}
\caption{Action de $\tilde{n}$}
\label{rev plan projectif}

\end{figure}

\begin{figure}[Ht]
\centering
\psfrag{g3}{$\g'_3$}
\psfrag{g4}{$\g'_4$}
\psfrag{f_1}{$b_{1}$}
\psfrag{f2}{$b_{2}$}

\includegraphics[width=14cm,keepaspectratio=true]{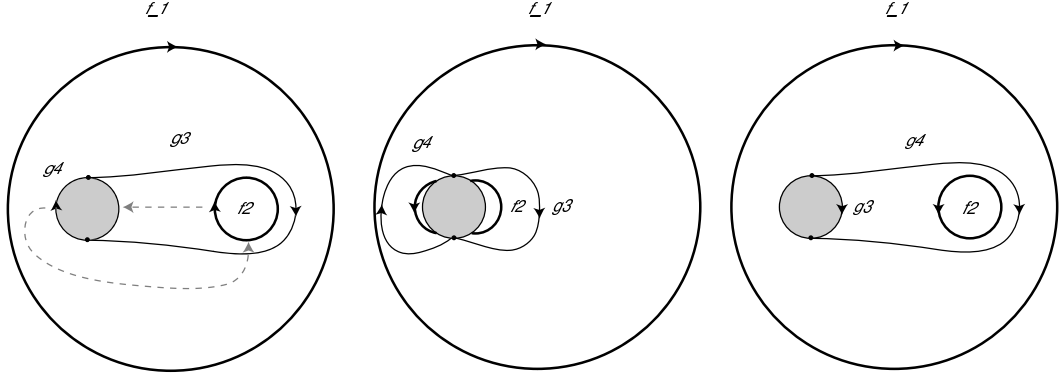}
\caption{Glissement de bord}
\label{glissement_bord}

\end{figure}

 Nous donnons deux représentations de $n$. La première (figure~\ref{rev plan projectif})
est l'action de $\tilde{n}$, un relevé de $n$, sur le revêtement des orientations $\tilde{Z}$ de $Z$.
$\tilde{Z}$ s'obtient, par exemple, en découpant $\g'_3$ dans $Z$ et en collant deux copies du pantalon obtenu
avec un twist de 1/2. En figure~\ref{rev plan projectif}), $\tilde{Z}$ est vu plongé dans l'espace $\R^3$,
$\tilde{n}$ est alors la reflexion par rapport au plan $\mathcal{P}$ contenant les
points de Weierstrass. Nous laissons au lecteur le soin de voir que ceci a bien un sens, quitte à donner
suffisament de symétries au plongement de $\tilde{Z}$.
 La deuxième représente un glissement de bord (figure~\ref{glissement_bord}) effectué par le bord $b_2$ en
suivant le chemin pointillé. Nous voyons bien que les géodésiques
intérieures de $Z$ sont échangées, que le bord $b_1$ conserve son orientation, et que celle de $b_2$
est inversée.\medskip

 Prenons $l_{\g'_3}$ comme coordonnée de $\t_{1,2}^-$, $\g'_4$ est l'autre géodésique fermée simple de $Z$.
L'action de $n$ s'exprime simplement en disant qu'au point de coordonnée $(l_{\g'_3})$ elle associe le point
de coordonnée $(l_{\g'_4})$~; $n$ a donc un unique point fixe, $Z(P)$, déterminé par $l_{\g'_3}=l_{\g'_4}$.
En tant que point fixe isolé de $n$, $Z(P)$ est eutactique. Il n'y a pas d'autre point eutactique car,
quitte à appliquer $n$, ce point aurait $\g'_3$ pour unique systole, or le gradient de $\g'_3$ est non nul.
En fait, nous voyons que $l_{\g'_4}$ est une fonction strictement décroissante de $l_{\g'_3}$.
Trivialement, $Z(P)$ est le seul point parfait de $\t_{1,2}^-$ et aussi le seul point extrême.\medskip

 Calculons la systole de $Z(P)$, comme $l_{\g'_3}=l_{\g'_4}$ nous pouvons considérer que $\bar{Z}(P)$ est un
élément du plan $\lbrace \theta_1=0\rbrace$, et de ce fait modélisée sur la figure~\ref{pentagone}.
Toute l'astuce du calcul est d'utiliser la relation trigonométrique entre une géodésique d'une surface$_{-1}$
et sa duale (rem.\ref{rem dualite}). Nous avons
\begin{eqnarray*}
\cosh(l_3/2) & = & \cosh(l_2/2)\cosh(l_1/2) \\
           & = & \sinh(l_{\g'_1}/2)\sinh(l_X/2) \sinh(l_{\g'_2}/2)\sinh(l_X/2), \\
\end{eqnarray*}
puis
\begin{eqnarray*}
\sinh^2(l_{\g'_3}/2) & = & \frac{\cosh^2(l_3/2)}{\sinh^2(l_X/2)} \\
                  & = & \sinh^2(l_{\g'_1}/2) \sinh^2(l_{\g'_2}/2) \sinh^2(l_X/2), \\
\end{eqnarray*}
or
\begin{eqnarray*}
\cosh(l_X/2) & = & \coth(l_{\g'_1}/2)\coth(l_{\g'_2}/2), \\
\end{eqnarray*}
ainsi
\begin{eqnarray*}
\sinh^2(l_{\g'_3}/2) & = & \sinh^2(l_{\g'_1}/2) \sinh^2(l_{\g'_2}/2) [\coth^2(l_{\g'_1}/2)\coth^2(l_{\g'_2}/2) - 1] \\
                  & = & \cosh^2(l_{\g'_1}/2)\cosh^2(l_{\g'_2}/2) - \sinh(l_{\g'_1}/2)\sinh(l_{\g'_2}/2) \\
                  & = & \cosh(\frac{l_{\g'_1}-l_{\g'_2}}{2})\cosh(\frac{l_{\g'_1}+l_{\g'_2}}{2}) \\
\end{eqnarray*}
ou encore
\begin{eqnarray*}
\cosh^2(l_{\g'_3}) & = & \cosh^2(l_{\g'_1}/2) + \cosh^2(l_{\g'_2}/2). \\
\end{eqnarray*}
Finalement nous avons démontré le théorème suivant~:

\begin{thm}
Modulo l'action du groupe modulaire, le maximum de la fonction systole sur $\t_{1,2}^-$ est atteint au point $Z(P)$, sa valeur est donnée par
$$\cosh(sys(Z(P)))=\cosh(l(b_1)/2)+\cosh(l(b_2)/2)+1.$$
\end{thm}

\section{Autres invariants métriques}
Nous donnons ici d'autres inégalités optimales valables pour la systole orientable, la systole non orientable, la 2-systole et la 3-systole des surfaces$_{-1}$. Tout consiste à exploiter convenablement des décompositions du domaine $D$.\par
\begin{figure}[Ht]
\psfrag{t1}{$\theta_1$}
\psfrag{z}{$z$}
\psfrag{x}{$x$}
\psfrag{C1}{$E1$}
\psfrag{C2}{$E2$}
\psfrag{C3}{$E3$}
\psfrag{C4}{$E4$}
\psfrag{C5}{$E5$}
\psfrag{C6}{$E6$}

\psfrag{lxl1}{$\lbrace l_1=l(\g_X) \rbrace$}
\psfrag{lxd1}{$\lbrace l(\g_X)=l(\g'_1) \rbrace$}
\psfrag{l1d1}{$\lbrace l_1=l(\g'_1) \rbrace$}
\psfrag{l1l2}{$\lbrace l_1=l_2 \rbrace$}

\includegraphics[width=12cm,keepaspectratio=true]{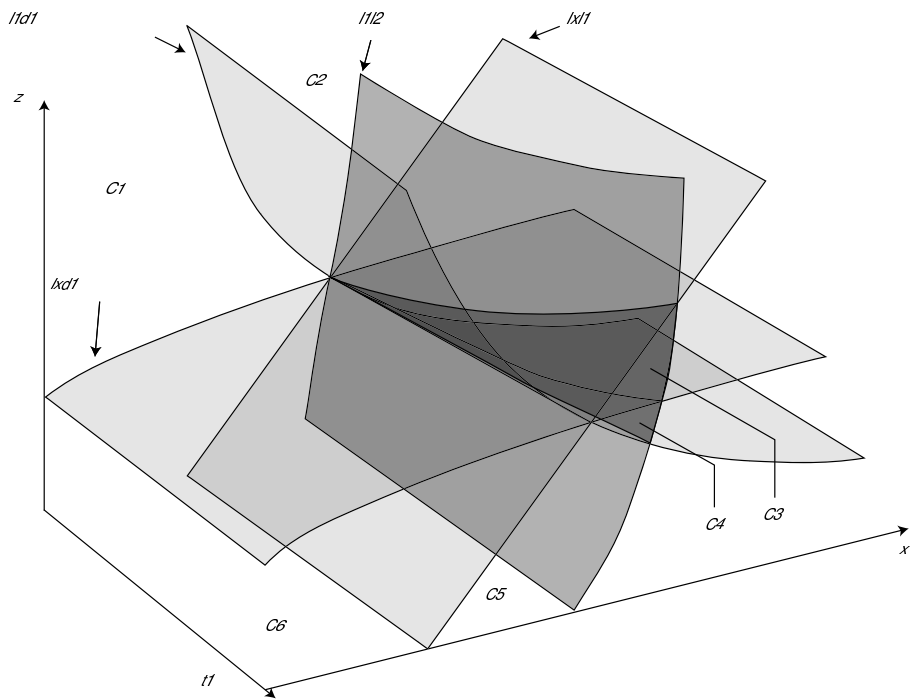}
\caption{Une deuxième décomposition cellulaire}
\label{decomposition2}
\end{figure}

 Nous pouvons par exemple choisir une décomposition avec six cellules de dimension 3~: 
\begin{itemize}
\item $E_1=\lbrace [X,\f]\in D\ |\ l_{\g_X}>l_{\g'_1}>l_{\g_1} \rbrace$, 
\item $E_2=\lbrace [X,\f]\in D\ |\ l_{\g_X}>l_{\g_1}>l_{\g'_1} \rbrace$,
\item $E_3=\lbrace [X,\f]\in D\ |\ l_{\g_1}>l_{\g_X}>l_{\g'_1} \rbrace$,
\item $E_4=\lbrace [X,\f]\in D\ |\ l_{\g_1}>l_{\g'_1}>l_{\g_X} \rbrace$, 
\item $E_5=\lbrace [X,\f]\in D\ |\ l_{\g'_1}>l_{\g_1}>l_{\g_X} \rbrace$, 
\item $E_6=\lbrace [X,\f]\in D\ |\ l_{\g'_1}>l_{\g_X}>l_{\g_1} \rbrace$.
\end{itemize}
cela correspond à diviser en deux chacune des cellules $C_i$ ($i=1,2,3$) en prolongeant naturellement les faces $F_{12}$, $F_{23}$ et $F_{13}$ (voir figure~\ref{decomposition2}).

\subsection{La systole orientable}
La systole orientable d'une surface$_{-1}$ est la plus petite longueur d'une géodésique orientable, qui éventuellement peut être carré d'une géodésique non orientable.
Dans $D$ il y a 3 candidats possibles pour la systole orientable~: $\g_1$, ${\g_X}^2$ ou ${\g'_1}^2$. Nous décomposons facilement $D$ en trois zones correpondant aux lieus où les différents candidats sont effectivement réalisés comme la systole orientable. Ces trois zones se déterminent sans difficultés puisque~:
$$\begin{array}{clrcl}
l_{\g_1}=2l_{\g_X}  & \Leftrightarrow & z       & = & \frac{\sqrt{x}+1}{2}  \\
l_{\g'_1}=l_{\g_X}  & \Leftrightarrow & (z-1)^2 & = & x  \\
l_{\g_1}=2l_{\g'_1} & \Leftrightarrow & z       & = & 1+\frac{2x}{\sqrt{x}-1} .\\
\end{array}$$
avec $x=\cosh^2(l_1/2)$ et $z=\cosh^2(l_X/2)$.
Nous regardons le maximum de la systole orientable dans chacune de ces zones, il est toujours atteint sur la droite des surfaces équilatérales $\lbrace l_1=l_2, \theta_1=1/2 \rbrace$. Sur cette droite il y a donc deux points intéressants, comme $\g_1$ est systole orientable de chacun d'eux nous voyons que le point réalisant le maximum global de la systole orientable sur $D$ est celui admetant $\g_1$, ${\g'_1}^2$, $\g_2$, ${\g'_2}^2$, $\g_3$, ${\g'_3}^2$ comme systoles orientables. Pour calculer la valeur de la systole orientable en ce point nous utilisons les égalités
$$\begin{array}{llll}
\cosh(h/2) &=& \coth(l_1/4)\coth(l_1/2) & \\
\cosh(l_1/2) &=& \cosh(h/2)\cosh(l_1/4) & \coth(l_1/2) \cosh^2(l_1/4)/\sinh(l_1/4)
\end{array}$$
où $h$ est la perpendiculaire commune aux bords associés à $\g_1$ dans le pantalon issu de la découpe de $\g_1$ dans $\T_X$. Partant de
$$\sinh(l_1/2)\sinh(l_1/4)=\cosh^2(l_1/4)$$
nous arrivons à l'équation
$$2X^2-5X+1=0$$
où $X=\cosh(l_1/2)$. Nous trouvons donc $X=(5+\sqrt{17})/4\cong 2,28$. Pour mémoire le cosinus hyperbolique de la demi-systole orientable de $X(H)$ vaut $(3+\sqrt{17})/4$. Nous pouvons énoncé le résultat suivant~:

\begin{pro}
Modulo l'action du groupe modulaire, la systole orientable atteint son maximum sur $\t_3^-$ au point de twist 1/2 tel que les géodésiques $\g_1$, ${\g'_1}^2$, $\g_2$, ${\g'_2}^2$, $\g_3$, ${\g'_3}^2$ aient mêmes longueurs. En ce point nous avons
$$\cosh(l_1/2)=(5+\sqrt{17})/2.$$
\end{pro}

\subsection{La systole non orientable}\label{systole non orientable}

Dans un premier temps, regardons seulement les géodésiques non orientables distinctes de l'ovale. Il est facile de voir que le minimum des fonctions longueur associées à ces géodésiques n'est pas borné sur $\t_3^-$. En effet, dans le domaine $D$ ce minimum est égal à $l_{\g'_1}$, et nous avons $\sinh(l_{\g'_1}/2)=\cosh(l_1/2)/\sinh(l_X/2)\geq 1/\sinh(l_X/2)$ qui tend vers $+\infty$ lorsque $l_X$ tend vers 0.\par
 Si maintenant nous considérons toutes les géodésiques non orientables, y compris l'ovale. Alors la systole non orientable (la plus petite longueur d'une géodésique non orientable) est égale sur $D$ au minimum des fonctions $l_X$ et $l_{\g'_1}$. Sur le sous-ensemble $\lbrace l_X \leq l_{\g'_1} \rbrace$ de $D$ le maximum de la systole non orientable, c'est-à-dire $l_X$, est atteint au point $X(H)$. Sur le sous-ensemble $\lbrace l_X \leq l_{\g'_1} \rbrace$ de $D$, à $l_X$ fixé le maximum de $l_{\g'_1}$ est atteint là où $l_1$ est maximum. Il y a deux cas à distinguer~: dans le premier cas le plan horizontal déterminé par la valeur de $l_X$ intersecte l'ensemble $\lbrace l_X=l_{\g'_1} \rbrace$ et $l_{\g'_1}\leq l_X \leq l_X(X(H))$~; dans le deuxième cas le maximum de la systole non orientable sur le plan horizontal est atteint à l'intersection des plans $\lbrace \theta_1=1/2 \rbrace$ et $\lbrace l_1=l_2 \rbrace$ (c'est la droite des surfaces équilatérales). Cependant à twist 1/2 nous avons,
$$\left\lbrace \begin{array}{lll}
\sinh(h/2) & = & \cosh(l_X/2)/\sinh(l_1/2) \\
\cosh(l_2/2) & = & \cosh(l_1/4)\cosh(h/2)
\end{array} ,\right. $$
 en notant $h$ la perpendiculaire commune aux bords associés à $\g_1$ dans le pantalon provenant de la découpe de $\g_1$ dans $\T_X$. Nous en déduisons
$$\cosh^2(l_2/2)=\frac{\cosh^2(l_X/2)+\cosh^2(l_1/2)-1}{2(\cosh(l_1/2)-1)},$$
et comme $l_1=l_2$ il vient
$$\sinh^2(l_X/2)=\cosh^2(l_1/2)(2\cosh(l_1/2)-3),$$
d'où
$$\sinh^2(l_{\g'_1}/2) = \frac{\cosh^2(l_1/2)}{\sinh^2(l_X/2)} = \frac{1}{2\cosh(l_1/2)-3}.$$
Le minimum de $\cosh(l_1/2)$ sur l'ensemble $\lbrace \theta_1=1/2, l_1=l_2, l_X \leq l_{\g'_1} \rbrace$ est atteint au point $X(H)$. Nous sommes arrivés au résultat suivant~:

\begin{pro}
Modulo l'action du groupe modulaire, la systole non orientable atteint son maximum sur $\t_3^-$ au point $X(H)$.
Cette surface possède 4 systoles non orientables~: $\g'_1$, $\g'_2$, $\g'_3$ et $\g_X$~; leur longueur est donnée par $$\cosh(l_1(X(H)))=(5+\sqrt{17})/2.$$
\end{pro}

\subsection{La 2-systole} 
La $k$-systole fut introduite par P. Scmutz Schaller dans \cite{schmutz}, nous renvoyons à cet article pour les définitions relatives à cet objet.

\begin{pro}
Modulo l'action du groupe modulaire, la 2-systole atteint son maximum global sur $\t_3^-$ en la surface équilatérale donnée par les conditions $l_1=l_2=l_3=l_{\g'_1}=l_{\g'_2}=l_{\g'_3}$. En ce point la systole est $\g_X$ et la 2-systole est réalisée par les 6 géodésiques $\g_i, \g'_i$ ($i=1,2,3$). Soit $X=\cosh(l_1/2)$, alors $2X^3-3X-2X+2=0$, en valeur approchée nous avons $\cosh(l_1/2)\cong 1,74$.
\end{pro}
\medskip

\begin{dem} 
 Un 2-système consiste en une géodésique orientable et sa duale, ou en une géodésique orientable et l'ovale, ou en deux géodésiques non orientables (distinctes de l'ovale) disjointes. Pour une surface fixée, si $l_X, l_1\leq l'_1$ alors $\lbrace l_X, l_1 \rbrace$ forme un 2-système de longueur minimale.\par
 Dans l'ensemble des points de $D$ vérifiant $l_X\leq l_1\leq l_{\g'_1}$ la 2-systole est donnée par $l_1$. Sur cet ensemble elle atteint son maximum au point de twist 1/2 tel que $l_1=l_{\g'_1}=l_2=l_{\g'_2}$, nous noterons $M_2$ ce point. Remarquons que la systole de $M_2$ est $\g_X$. La 2-systole de $M_2$, c'est-à-dire $l_1(M_2)$, vérifie l'équation $2X^3-3X-2X+2=0$ avec $X=\cosh(l_1(M_2)/2)$.
Nous trouvons $\cosh(l_1(M_2)/2)\cong 1,74$.\par
 Dans l'ensemble des points vérifiant $l_1\leq l_X\leq l_{\g'_1}$, la 2-systole est donnée par $l_X$. Son maximum est atteint par exemple au point $X(P)$, et vaut 
 $$\cosh(l_X(X(P))/2)=(1+\sqrt{5})/2\cong 1,62.$$\par
 Dans l'ensemble des point vérifiant $l_1\leq l_{\g'_1} \leq l_X$, la 2-systole est donnée par $l_{\g'_1}$. Son maximum est atteint lui aussi (par exemple) au point $X(P)$, et a donc même valeur que précédemment.\par 
 Enfin, dans tous les autres cas, c'est-à-dire lorsque $l_{\g'_1}\leq l_1$, la 2-systole est la géodésique de plus petite longueur parmi celles disjointes de $\g'_1$. Clairement, le cosinus hyperbolique de la 2-systole sur cet ensemble est donné par $\cosh(l_{\g'_1}/2)+1$ (formule de la systole de la bouteille de Klein à bord). Donc le maximum de la 2-systole sera atteint là où $l_{\g'_1}$ atteint son maximum. A $l_X$ fixée, $l_{\g'_1}$ atteint son maximum là où $l_1$ est maximum, soit à l'intersection des plans $\lbrace \theta_1=1/2 \rbrace$ et $\lbrace l_1=l_2 \rbrace$ (la droite des surfaces équilatérales). Nous reproduisons à l'identique le raisonnement du paragraphe précédent, nous voyons que le minimum de $\cosh(l_1/2)$ sur l'ensemble $\lbrace \theta_1=1/2, l_1=l_2, l_{\g'_1}\leq l_1 \rbrace$ est atteint au point $M_2$.\par
 Nous concluons que, modulo l'action de $Mod_3^-$, le maximum de la 2-systole sur $\t_3^-$ est atteint au point $M_2$.\end{dem}\hfill$\square$


\subsection{La 3-systole} 
Un 3-système d'une surface$_{-1}$ est constitué de 3 géodésiques disjointes, il s'en suit que ces géodésiques sont nécessairement non orientables et distinctes de l'ovale. Nous avons vu au \textsection~\ref{systole non orientable} que le minimum des géodésiques non orientables distinctes de l'ovale n'est pas borné sur $\t_3^-$, par conséquent la 3-systole elle non plus n'est pas bornée sur $\t_3^-$.

\newpage
\nocite{*}
\bibliographystyle{alpha}
\bibliography{systole_carac_-1}
\hspace{1cm}

\noindent Matthieu Raphaël GENDULPHE\\
Laboratoire Bordelais d'Analyse et Géométrie\\ 
Université Bordeaux1\\
351, cours de la Libération\\
33 405 Talence, FRANCE\\
matthieu.gendulphe@math.u-bordeaux1.fr
\end{document}